\documentclass[10 pt, a4paper]{article}

\usepackage{graphicx}
\usepackage[utf8]{inputenc}  
\usepackage{amssymb}
\usepackage{amsmath}
\usepackage{hyperref}
\usepackage{xfrac}
\usepackage{times}
\usepackage{url}
\usepackage{adjustbox}
\usepackage{boldline}
\usepackage[T1]{fontenc}
\usepackage[dvips]{epsfig}
\usepackage{epstopdf}
\usepackage{color}
\DeclareGraphicsExtensions{.ps,.eps,.pdf}
\allowdisplaybreaks

\begin{document}
\pdfoutput=1
\date{}

\title{Identification of Induction Motors with Smart Circuit Breakers}

\author{Lorenzo~Fagiano,
	Marco Lauricella,
	Daniele Angelosante,
	and Enrico Ragaini%
	\thanks{This manuscript is the pre-print of a paper submitted for possible publication on the IEEE Transactions on Control Systems Technology and it is subject to IEEE Copyright. If accepted, the copy of record will be available at IEEEXplore library: \url{http://ieeexplore.ieee.org/}.}
	\thanks{L. Fagiano and M. Lauricella are with the Dipartimento di Elettronica, Informazione e Bioingegneria, Politecnico di Milano. D. Angelosante is with ABB Switzerland, Corporate Research. E. Ragaini is with ABB Italia SpA. E-mail addresses: \{lorenzo.fagiano $|$ marco.lauricella\}@polimi.it; daniele.angelosante@ch.abb.com; enrico.ragaini@it.abb.com.}
	\thanks{Corresponding author: Lorenzo Fagiano.}}

\maketitle

\begin{abstract}
The problem of estimating the parameters of induction motor models is considered, using the data measured by a circuit breaker equipped with industrial sensors. The measured data pertain to direct-on-line motor startups, during which the breaker acquires three-phase stator voltage and current derivative. This setup is novel with respect to previous contributions in the literature, where voltage and current (and possibly also rotor speed) are considered. The collected data are used to formulate a parameter identification problem, where the cost function penalizes the discrepancy between simulated and measured derivatives of the stator currents. The resulting nonlinear program is solved via numerical optimization, and a number of algorithmic improvements with respect to the literature are proposed. In order to evaluate the goodness of the obtained results, an experimental rig has been built, where the motor's voltages and currents are simultaneously acquired also by accurate sensors, and the corresponding identification results are compared with those obtained with the circuit breaker. The presented experimental results indicate that the considered industrial circuit breaker is able to provide data with high-enough quality to carry out model-based nonlinear identification of induction machines. The identified models can then be used for several further applications within a smart grid scenario.
\end{abstract}


\section{Introduction}
\label{S:intro}
The proliferation of sensing and computing devices and communication capabilities in power networks is at the basis of the smart grid paradigm \cite{Fang2012}, where bi-directional flows of electricity and information are exploited to improve and automate grid operation and to enable the use of distributed electricity generation. Among the different functions that a smart grid shall accomplish, self-monitoring, self-healing, and  advanced load protection and monitoring are crucially important to limit operation and maintenance costs as the grid complexity increases. However, to realize these functions in a capillary way requires the installation and connection of a large number of sensing devices, thus adding more complexity and costs.\\ 
Circuit breakers represent an ideal candidate to alleviate this problem. Installed in millions across the power grid at all voltage levels, these devices are designed to last tens of years and efficiently protect portions of the grid, groups of loads and individual loads from the effects of short circuits and other types of faults. Traditionally built as passive and isolated devices, particularly in low voltage installations, circuit breakers can provide a distributed network of sensors and actuators if equipped with sensing, computing and communication capabilities. Since the breakers are already connected to the grid, there is no need to install a separate link to power the sensors and on-board processors. These smart circuit breakers can then accomplish additional functionalities with respect to the classical protection one. An example is the ABB Emax2$^\circledR$ circuit breaker, which can also operate as power manager by selectively disconnecting downstream loads in order to avoid an excessive power consumption in each billing period \cite{ABB2014}.

In addition to energy management, another function of interest that could be carried out by smart breakers is the identification of suitable models of the connected loads, using the measurements of their electric signature. Then, the identified models can be used for purposes such as load detection and monitoring and discrimination between loads with high inrush currents and faults in the network, as proposed in the patent \cite{Ragaini2017} related to the research presented here. 
In this paper, we explore this functionality by considering an industrial scenario, where the most common type of load is represented by electric motors, accounting for about 69\% of the whole electricity consumption of the industry sector \cite{Waide2011}. In particular, three-phase asynchronous alternating current (AC) induction motors with direct on-line (DOL) connection are most frequently used and they are considered here.\\
In the described context, the most relevant question to be answered is whether the data collected by industrial current and voltage sensors, installed in commercial circuit breakers, are good enough to carry out the parametric identification of an induction motor's model from its electrical signature. The main contribution of our work is to show that indeed this is possible already now. This claim is supported by extensive experimental tests where we compare the results obtained with industrial sensors with those obtained with highly accurate laboratory sensors.\\
This outcome indicates that the currently marketed circuit breakers are already equipped with sensors and data acquisition systems suitable to implement the described model identification function. The identification algorithm itself can then run either locally, using the circuit breaker control unit, or remotely, for example through a cloud service to which the smart breaker sends the measured data and receives back the estimated parameters' values. The problem of identifying the lumped parameters of a motor's model from the measured three-phase currents and voltages has been addressed in the literature by several contributions, using e.g. Recursive Least-Squares algorithms \cite{Koubaa2004}, Genetic Algorithms \cite{Alonge1998}, Extended Kalman Filters \cite{Iwasaki1989} or Total Least Squares approaches plus neurons \cite{Cirrincione2003new}. In this work, the parameter identification problem is cast into a nonlinear least squares estimation, where a batch of data collected during the motor direct on-line start-up transient is compared with the simulated quantities, obtained by integrating the model from known initial conditions and applying the measured stator voltage. The resulting optimization problem is a nonlinear program, which is commonly treated with Newton-type iterative optimization algorithms. In the literature, this approach has been already considered, with different choices of measured variables: the most common setup is to assume the availability of measurements of the stator's currents and voltages and of the rotor's angular speed, see e.g. \cite{Wang2005nl}, \cite{Cirrincione2005cls}, \cite{Proco1999}. In fewer cases, only currents and voltages are used, like in \cite{Shaw1999}, \cite{Velez2001cdc}. In all these works, aspects like  sensitivity of the estimation procedure to initialization,  stability of the estimated parameters with respect to the chosen sampling frequency, and efficiency of the employed optimization routine are not fully discussed. Yet, these are crucial issues from the point of view of control system technology implementation. As additional contributions, in this paper we present several results concerning these points, in particular we show how in this problem a Constrained Gauss-Newton optimization method \cite{NoWr06} combined with analytical computation of sensitivities leads to a rather efficient solution approach, and how using a trapezoidal numerical integration yields good stability of the estimation algorithm also with rather low sampling frequency, poor initialization, and model over-parametrization, especially when compared with the (most commonly used) forward Euler integration method. These contributions are also novel with respect to our recent work \cite{Angelosante2017}, in which we considered forward Euler integration and we did no attempt to improve the efficiency of the identification routine. Finally, a further novelty introduced in our research is the direct use, in the identification routine, of data obtained with Rogowski coil sensors (installed in the considered commercial circuit breakers), which measure the current derivative instead of the current itself (as considered so far in the literature). 

The paper is organized as follows. Section \ref{S:probl_setup_descr} introduces the experimental setup that we built to carry out our tests and the formulation of the parameter identification problem. Section \ref{S:probl_model_data} presents the chosen induction machine model and describes the available measurements. Section \ref{S:nonlinid} describes the considered nonlinear identification approach and the proposed improvements. Experimental results are discussed in Section \ref{S:results}, and conclusions and future developments in Section \ref{S:conclusions}.

\begin{figure}[thpb]
	\centering
	\includegraphics[width=0.65\columnwidth]{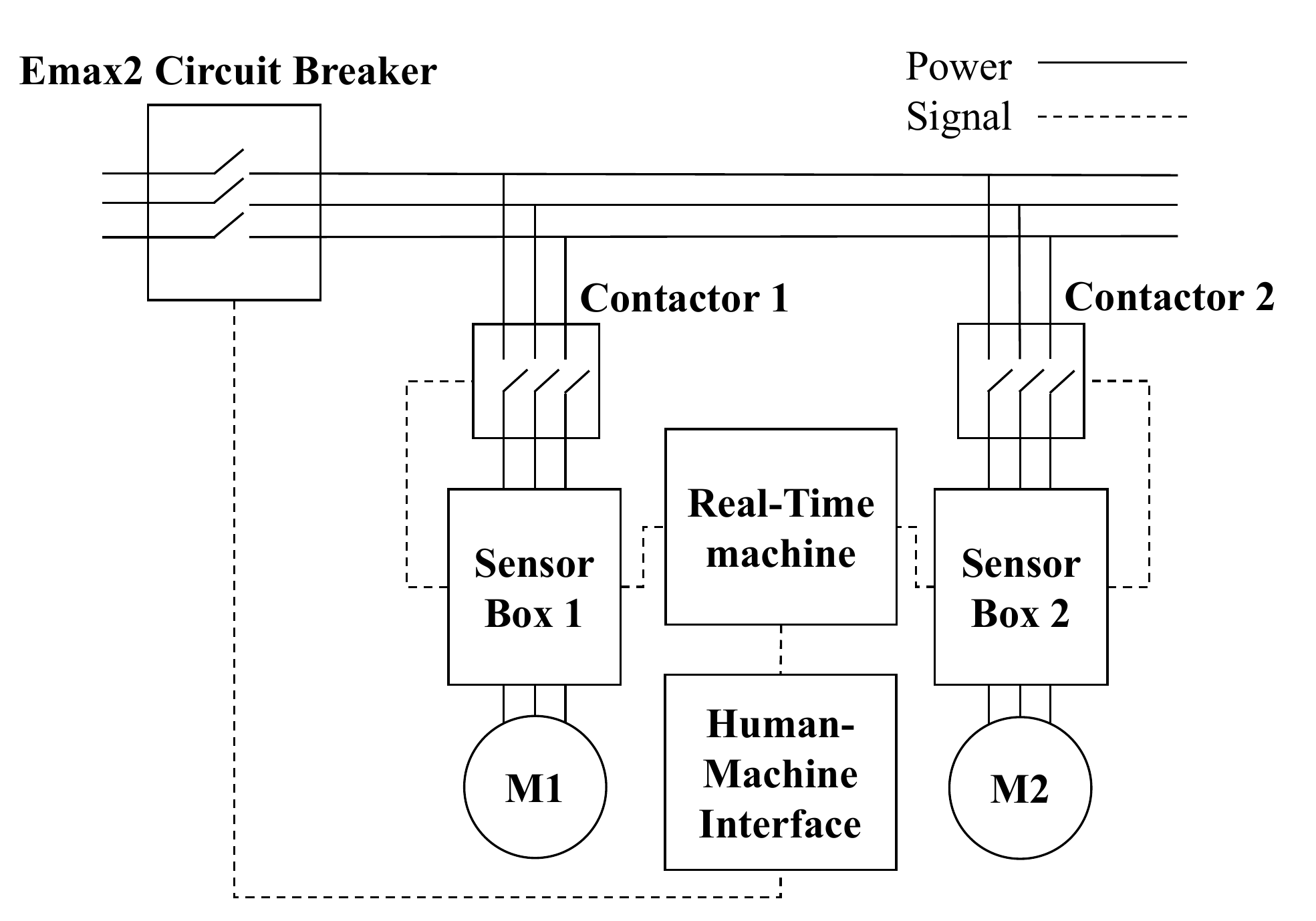}
	\caption{Layout of the employed experimental setup. An Emax2$^\circledR$ circuit breaker is installed upstream of two electric motors, which are connected through remotely piloted contactors. The breaker acquires measurements of the three-phase voltage and of the aggregate three-phase current derivatives. Two sensors boxes, one for each motor, acquire high-quality voltage and current measurements. The motors' switching is governed by a real-time machine, operated by the testing personnel. The open/close signals sent by the real-time machine to each motor are relayed by the sensor boxes to the contactors. The real-time machine also logs the data collected by the breaker and by the sensor boxes. 220-VAC supply circuits for the real-time machine and the sensor boxes are not shown in this figure.}
	\label{f:exp_layout}
\end{figure} 
\begin{figure}[thpb]
	\centering
	\includegraphics[width=0.65\columnwidth]{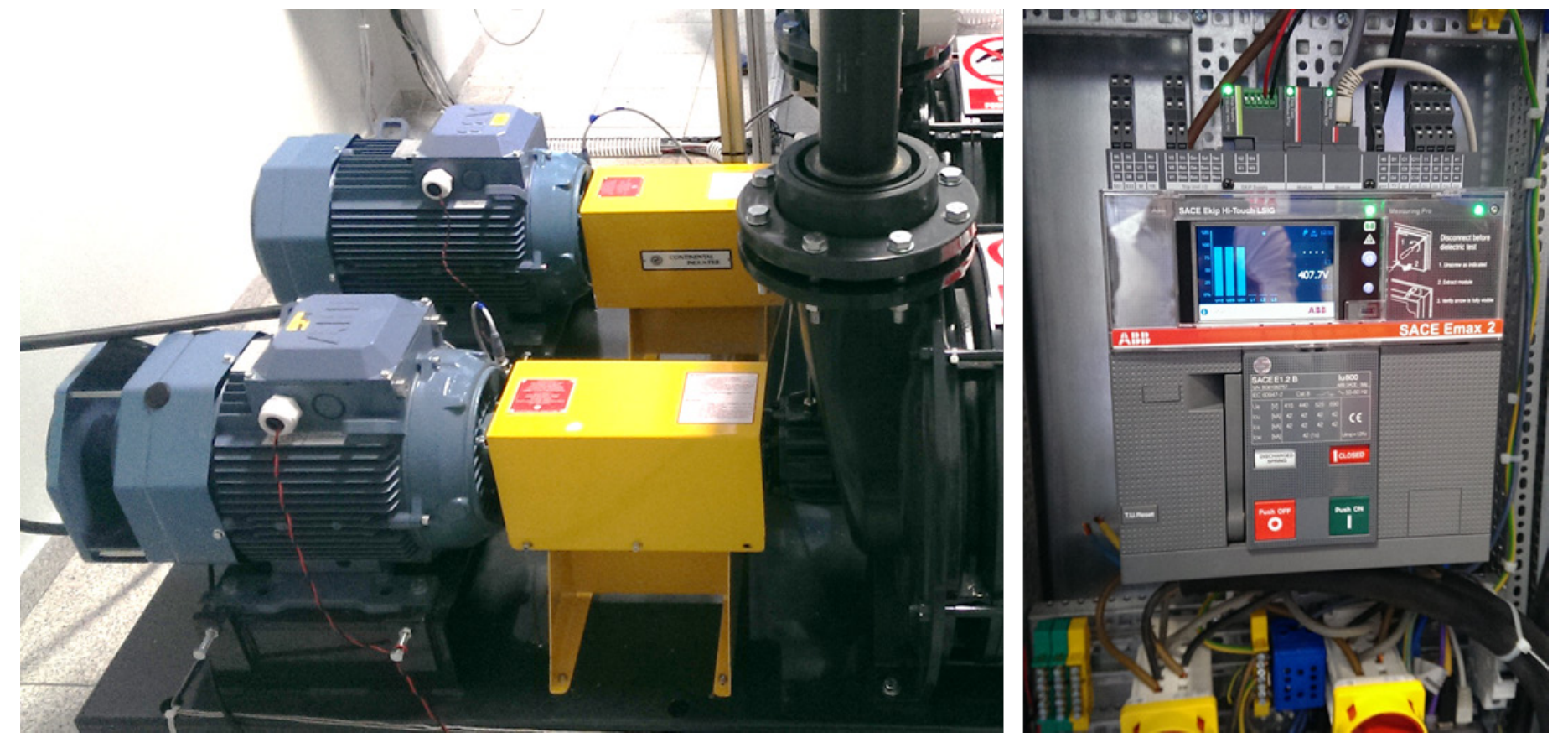}
	\caption{Pictures of the experimental setup. Left: induction motors employed for the tests. Right: Emax2$^\circledR$ breaker installed in the electric cabinet of the testing laboratory.}
	\label{f:exp_picture}
\end{figure} 

\section{System description and problem formulation}
\label{S:probl_setup_descr}
The layout of the considered experimental setup is shown in Fig. \ref{f:exp_layout}. A picture of the setup, realized at ABB Corporate Research in Poland, is shown in Fig. \ref{f:exp_picture}. The network voltage is 380 VAC phase-to-phase, 50 Hz.  Referring to Fig. \ref{f:exp_layout}, the system includes the following elements:
\begin{itemize}
	\item An ABB Emax2$^\circledR$ circuit breaker with 800 A (rms) of nominal current, which measures the three phase voltage (phase-to-phase measurements) via a resistive divider, and the three-phase current derivatives via Rogowski coils;
	\item Two induction motors, $M1$ and $M2$, and two contactors (ABB AF38 series) to connect each motor to the 3-phase line. Motor $M1$ is Y-connected, while $M2$ is $\Delta$-connected.
	\item Two sensor boxes, built at ABB Corporate Research Switzerland, equipped with three-phase voltage and current sensors based on Hall-effect transducers, and with a relay to send open/close control signals to the contactors;
	\item A Real-Time machine, which logs the data acquired by the circuit breakers and the sensor boxes, and sends the open/close commands to the contactors  via the sensor boxes. The Real-Time machine is operated by the testing personnel via a Human-Machine Interface, to carry out the desired testing sequences.
\end{itemize} 
The main features of the employed motors, sensors, and data acquisition systems are presented in Table \ref{t:sensor_specs}. The sensor boxes feature highly accurate transducers and the corresponding measured data is used, in our research, as ``ground truth'' to evaluate the performance that can be achieved with the data collected by the circuit breaker, which is the object of study.  

\begin{table}[thpb]
	\caption{Main paramaters of the employed motors, sensors, and data acquisition systems.}
	\label{t:sensor_specs}
	\centering
	\begin{adjustbox}{max width=1\columnwidth}
		\begin{tabular}{|l|l|}
			\hline
			\multicolumn{2}{|c|}{\textbf{Motor $M1$}}\\\hline
			Model& ABB M3AA160MLB2\\\hline
			Rated power& 15 kW\\\hline
			N. of poles& 2\\\hline
			Connection& Star\\\hline\hline
			\multicolumn{2}{|c|}{\textbf{Motor $M2$}}\\\hline
			Model& ABB M3AA160MLC23GAA\\\hline
			Rated power& 18 kW\\\hline
			N. of poles& 2\\\hline
			Connection& Delta\\\hline\hline
			\multicolumn{2}{|c|}{\textbf{Emax2$^\circledR$ Circuit Breaker}}\\\hline
			\multicolumn{2}{|c|}{3-phase voltage measurements}\\\hline
			Transducer type& Resistive divider\\\hline
			Range& 690 Vrms, phase-to-phase\\\hline
			Bandwidth& 1 kHz\\\hline
			\multicolumn{2}{|c|}{3-phase current derivative measurements}\\\hline
		Transducer type& Rogowski coil\\\hline
		Range& 65$\,$10$^6\,$As$^{-1}$\\\hline
		Bandwidth& 1 kHz\\\hline
		\multicolumn{2}{|c|}{Data-acquisition system}\\\hline
		N. of bits&12\\\hline
		Sampling frequencies& 1.2, 2.4, 4.8, 9.6 kHz\\\hline\hline
		\multicolumn{2}{|c|}{\textbf{Sensor boxes}}\\\hline
		\multicolumn{2}{|c|}{3-phase voltage measurements}\\\hline
		Transducer type& Resistor + Hall effect, compensated (LEM LV 25-P)\\\hline
		Range& 400 Vrms, phase-to-phase or phase-to-neutral selectable\\\hline
		Bandwidth& 100 kHz\\\hline
		\multicolumn{2}{|c|}{3-phase current measurements}\\\hline
		Transducer type& Hall effect, compensated (LEM LF 205-S)\\\hline
		Range& $\pm$420 A\\\hline
		Bandwidth& 100 kHz\\\hline
		\multicolumn{2}{|c|}{Data-acquisition system}\\\hline
		N. of bits& 16\\\hline
		Sampling frequency& 5 kHz\\\hline
\end{tabular}
	\end{adjustbox}
\end{table}

The experimental tests considered in this work are direct-on-line motor startups, in which both contactors are initially open. Then, upon command by the test personnel, the Real-Time machine sends a triggering signal to one of the two contactors and acquires the electric signature of the corresponding motor, as measured both by its sensor box and by the smart breaker. This testing procedure is well-motivated by the possibility, in a real-world application, to carry out several motor startups in the commissioning phase of a new installation, in order to record the electric signature of each electric machine in a controlled way for the sake of parameter estimation. We built a setup with two induction machines in order to carry out also tests with simultaneous and/or partially overlapping startups, for the sake of studying motor detection and discrimination between different machines. These topics are subject of future research.

Given the data obtained in the start-up tests, the problem we address is to identify the parameters of a model of each motor, where the inputs are the stator voltages and the measured outputs are either the stator currents (for sensor box data) or their derivatives (for circuit breaker data). In particular, we seek the parameter values that minimize an output-error (i.e. simulation) performance criterion. The cost function is in fact based on the error between the measured outputs and those computed by simulating the model from known initial condition (standstill), by applying in open loop the measured input values. Once obtained, these parameters can be used for a number of further tasks such as real-time monitoring, filtering, and fault detection.\\ In the next section, we describe the chosen motor model and the measurements available for the parameter identification procedure.

\section{Induction machine model and experimental data-set}
\label{S:probl_model_data}
\subsection{Induction machine model}
\label{SS:model}
We resort to a rather standard dynamical model of three-phase induction motors, summarized here for the sake of completeness. The model is derived starting from the voltage balance equations of the stator and rotor phases, following the procedure described e.g. in \cite{Krause1986}. We assume linearity of the inductances and neglect the losses in the iron. \\
\noindent Let us define the following three-phase quantities:
\begin{itemize}
	\item stator voltage $\boldsymbol{v}_{abc,s}(t):=[v_{as}(t) \ v_{bs}(t) \ v_{cs}(t)]^T$
	\item stator current $\boldsymbol{i}_{abc,s}(t):=[i_{as}(t) \ i_{bs}(t) \ i_{cs}(t)]^T$
	\item rotor voltage $\boldsymbol{v}_{abc,r}(t):=[v_{ar}(t) \ v_{br}(t) \ v_{cr}(t)]^T$
	\item rotor current $\boldsymbol{i}_{abc,r}(t):=[i_{ar}(t) \ i_{br}(t) \ i_{cr}(t)]^T$
	\item  stator fluxes $\boldsymbol{\lambda}_{abc,s}(t):=[\lambda_{as}(t) \ \lambda_{bs}(t) \ \lambda_{cs}(t)]^T$
	\item rotor fluxes $\boldsymbol{\lambda}_{abc,r}(t):=[\lambda_{ar}(t) \ \lambda_{br}(t) \ \lambda_{cr}(t)]^T$
\end{itemize}
where $t$ is the continuous time variable, $a$, $b$, $c$ denote the motor phases, $s$ and $r$ denote the stator and the rotor quantities respectively and $\lambda$ is the total flux of a particular winding. 
It is customary to transform three-phase quantities into two-phase ones by resorting to a suitable change of variables which implies the choice of a common reference frame. Here, we select a fixed (i.e. non-rotating) frame which, under the assumption that the electric machine is balanced (i.e. $i_{as}(t)+i_{bs}(t)+i_{cs}(t)=0$ and $i_{ar}(t)+i_{br}(t)+i_{cr}(t)=0$) results in the following transformation matrix:
$$
\boldsymbol{M}:=\frac{2}{3} \begin{bmatrix}
1 & \cos{(-\frac{2}{3}\pi)} & \cos{(\frac{2}{3}\pi)} \\
0 & \sin{(-\frac{2}{3}\pi)} & \sin{(\frac{2}{3}\pi)} \\
0.5 & 0.5 & 0.5
\end{bmatrix}.
$$
The matrix $\boldsymbol{M}$, when applied to a generic vector in the form $\boldsymbol{s}_{abc}:=[s_a \ s_b \ s_c]^T$, pertaining to a balanced 3-phase frame, gives $\boldsymbol{s}_{dq0}=\boldsymbol{M}\boldsymbol{s}_{abc}$, where $\boldsymbol{s}_{dq0}=[s_d \ s_q \ 0]^T$, i.e. a vector of only two independent components, commonly referred to as the \textit{dq}-components. In the following, we denote as $\boldsymbol{s}_{dq}$ the 2-dimensional vectors in \textit{dq}-components where we dropped the zero component.\\ 
\noindent We further introduce the fluxes per time unit as $\psi_{ds}(t)=\omega_e \lambda_{ds}(t)$, $\psi_{qs}(t)=\omega_e \lambda_{qs}(t)$, $\psi_{dr}(t)=\omega_e \lambda_{dr}(t)$, $\psi_{qr}(t)=\omega_e \lambda_{qr}(t)$, where $\omega_e$ is the nominal electrical frequency of the AC network in rad/s. Moreover, the electrical torque $T_e(t)$ and the load torque $T_l(t)$ are modeled as:
\begin{equation}\label{E:elec_torque}
T_e(t)=\frac{3N_p}{4\omega_e}\left( \psi_{qr} i_{dr} - \psi_{dr} i_{qr} \right)
\end{equation}
\begin{equation}\label{E:load_torque}
T_l(t)=T_{l_0}+T_{l_1}\,w_r(t)
\end{equation}
where $\omega_r(t)$ is the rotor angular speed, $N_p$ is the number of poles of the motor, $J_r$ is the rotor moment of inertia, and $T_{l_0}$ and $T_{l_1}$ are, respectively, a constant load coefficient and a constant viscous friction coefficient. Note that this load model is over-parametrized with respect to our experimental setup, where the constant load is zero and only the linear viscous term is present. We show in our experimental results (Section \ref{SS:overparam}) that such over-parametrization is handled correctly by the estimation algorithm only for suitable choices of the sampling frequency and model discretization techniques.\\
\noindent In the model, we consider as input of the system the stator voltage in its \textit{dq}-representation, given by $\boldsymbol{u}(t):=\allowbreak [v_{ds}(t) \ v_{qs}(t)]^T \allowbreak \in\mathbb{R}^2$, and as state the vector $\boldsymbol{x}(t):=[\psi_{ds}(t) \ \psi_{qs}(t) \ \psi_{dr}(t) \ \psi_{qr}(t) \ \omega_r(t)]^T\in \mathbb{R}^5$. We are now in position to write the model equations pertaining to the state evolution in time:
\begin{equation}
\label{E:compactcontmatmodel}
\dot{\boldsymbol{x}}(t)=\boldsymbol{A}(\omega_r(t))\boldsymbol{x}(t)+\boldsymbol{B}\boldsymbol{u}(t)+\boldsymbol{\beta}(\boldsymbol{x}(t))
\end{equation}
where
\begin{equation}
\label{E:model_matrices}
\begin{array}{c}
\boldsymbol{A}(\omega_r(t))= \omega_e 
\begin{bmatrix}
\frac{R_s(X_{m}-X_{l})}{X_{l}^2} & 0 & \frac{R_sX_{m}}{X_{l}X_{l}} & 0 & 0 \\
0 & \frac{R_s(X_{m}-X_{l})}{X_{l}^2} & 0 & \frac{R_sX_{m}}{X_{l}X_{l}} & 0 \\
\frac{R_rX_{m}}{X_{l}X_{l}} & 0 & \frac{R_r(X_{m}-X_{l})}{X_{l}^2} & -\frac{\omega_r(t)}{\omega_e} & 0 \\
0 & \frac{R_rX_{m}}{X_{l}X_{l}} & \frac{\omega_r(t)}{\omega_e} & \frac{R_r(X_{m}-X_{l})}{X_{l}^2} & 0 \\
0 & 0 & 0 & 0 & 0
\end{bmatrix};\\
\boldsymbol{B}=
\begin{bmatrix}
\omega_e & 0 \\ 0 & \omega_e \\ 0 & 0 \\ 0 & 0 \\ 0 & 0
\end{bmatrix}; \quad
\boldsymbol{\beta}(\boldsymbol{x}(t))=
\begin{bmatrix}
0 \\ 0 \\ 0 \\ 0 \\ \frac{N_p}{2J_r}\left( T_e(t)-T_l(t) \right)
\end{bmatrix}.
\end{array}
\end{equation}
\noindent In \eqref{E:model_matrices}, $R_s$ is the stator resistance, $R_r$ is the rotor resistance, $X_{l}$ and $X_m$ are respectively the stator (and rotor) reactance and the magnetizing reactance at the nominal electric frequency. As usual, the notation $\dot{x}\doteq dx/dt$ denotes the time derivative. \\
Regarding the output equations, these depend on the measured output of the system, which can be either the stator current or its derivative (depending on the considered measuring equipment, see Sections \ref{S:probl_setup_descr} and \ref{SS:data-sets}). The output variables are again transformed in \textit{dq}-components. We indicate with $\boldsymbol{y}_{SB}(t)$ the output vector obtained with current measurements (i.e. the sensor boxes, see Section \ref{S:probl_setup_descr}) and with $\boldsymbol{y}_{CB}(t)$ the one given by current derivative measurements (i.e. the circuit breaker). Thus, in the first case we have:
\begin{equation}
\label{E:compactcontmatmodel_output_curr}
\boldsymbol{y}_{SB}(t)=\boldsymbol{C}\boldsymbol{x}(t),
\end{equation}
where
$$
\boldsymbol{C}= \frac{1}{X_{l}} \begin{bmatrix}
1-\frac{X_{m}}{X_{l}} & 0 & -\frac{X_{m}}{X_{l}} & 0 & 0 \\
0 & 1-\frac{X_{m}}{X_{l}} & 0 & -\frac{X_{m}}{X_{l}} & 0
\end{bmatrix}.
$$
If current derivative measurements are considered, we have:
\begin{equation}
\label{E:compactYDOTmatmodel}
\boldsymbol{y}_{CB}(t)=\boldsymbol{C}\dot{\boldsymbol{x}}(t).
\end{equation}
Equations \eqref{E:elec_torque}-\eqref{E:compactYDOTmatmodel} provide the continuous time model of the motor considered in this paper. The vector of parameters to be identified from experimental data is denoted with $\boldsymbol{p}=[R_s \ R_r \ X_l \ X_m \ J_r \ T_{l_0} \ T_{l_1}]^T,\,\boldsymbol{p}\in \mathbb{R}^7$, while the number of poles $N_p$ is assumed known since it is easily obtained from the motor nameplate or data-sheet.\\
We introduce next the considered available data-sets.
 
\subsection{Measured data sets}\label{SS:data-sets}
We indicate with  $t_s$ the sampling period (and with $f_s=1/t_s$ the sampling frequency) and with $N$ the total number of samples. For the sake of notational simplicity, we retain these two symbols irrespective of whether the considered data-sets are obtained with the sensor boxes or with the circuit breaker, since this aspect will be clear from the context. Moreover,  all the collected 3-phase measurements are transformed into \textit{dq}-components to be compatible with the employed model. Finally, as a general rule we indicate with $\tilde{\cdot}$ the measured (i.e. affected by noise) quantities. As regards the data collected by the sensor boxes, the measured voltage $\tilde{\boldsymbol{V}}_{SB}$ is given by:
$$
\tilde{\boldsymbol{V}}_{SB}= \begin{bmatrix}
\tilde{v}_{ds,SB}(t_s), \cdots, \tilde{v}_{ds,SB}(N\,t_s) \\
\tilde{v}_{qs,SB}(t_s), \cdots, \tilde{v}_{qs,SB}(N\,t_s)
\end{bmatrix},
$$
where $\tilde{v}_{ds,SB},\,\tilde{v}_{qs,SB}$ are the $dq$-components of the stator voltages acquired by the voltage transducers in the sensor boxes. Similarly, the measured current $\tilde{\boldsymbol{I}}$ is:
$$
 \tilde{\boldsymbol{I}}= \begin{bmatrix}
\tilde{i}_{ds}(t_s), \cdots, \tilde{i}_{ds}(Nt_s) \\
\tilde{i}_{qs}(t_s), \cdots, \tilde{i}_{qs}(Nt_s)
\end{bmatrix}.
$$
On the other hand, for the measurements collected by the smart circuit breaker we have:
$$
\tilde{\boldsymbol{V}}_{CB}= \begin{bmatrix}
\tilde{v}_{ds,CB}(t_s), \cdots, \tilde{v}_{ds,CB}(N\,t_s) \\
\tilde{v}_{qs,CB}(t_s), \cdots, \tilde{v}_{qs,CB}(N\,t_s)
\end{bmatrix},
$$
with $\tilde{v}_{ds,CB},\,\tilde{v}_{qs,CB}$ being the $dq$-components of the stator voltages acquired by the voltage transducers in the breaker, and
$$
\tilde{\dot{\boldsymbol{I}}}= \begin{bmatrix}
\tilde{\dot{i}}_{ds}(t_s), \cdots, \tilde{\dot{i}}_{ds}(Nt_s) \\
\tilde{\dot{i}}_{qs}(t_s), \cdots, \tilde{\dot{i}}_{qs}(Nt_s)
\end{bmatrix}.
$$
Notice that, differently from the voltage signals, for the measured current and current derivative signals  the acquisition device is unambiguous, since the sensor boxes measure only the current and the circuit breaker measures only its derivative. Therefore, for these we dropped the subscripts $_{CB}$ and $_{SB}$ for simplicity.
\section{Nonlinear identification procedure}
\label{S:nonlinid}
In this section we detail the estimation algorithm used to identify the induction motor parameters and the various proposed alternatives in terms of fitting criterion and computational aspects. As stated in Section \ref{S:intro}, the parameter identification problem is cast into a nonlinear least squares estimation, where a batch of data collected during the motor start-up transient is compared with the corresponding simulated quantities, obtained by integrating the model from known initial condition and applying the acquired stator voltage data. The resulting numerical optimization problem takes the general form:
\begin{subequations}
\label{E:opt_probl}
\begin{align}
&\hat{\boldsymbol{p}}=\underset{\boldsymbol{p}\in \mathcal{P}}{\arg \, \mathrm{min} \, tr} ( \boldsymbol{J}(\boldsymbol{p})) \label{E:opt_probl_cost}\\
&\text{subject to}\nonumber\\
&\text{discrete-time model equations}\label{E:opt_probl_eq}
\end{align}
\end{subequations}
where $tr(\cdot)$ indicates the trace of a matrix, $\boldsymbol{J}(\boldsymbol{p})$ is a square cost matrix, and $\mathcal{P}$ is a set of admissible parameters, defined e.g. by box constraints that account for reasonable upper and lower bounds on each component of $\boldsymbol{p}$. The set $\mathcal{P}$ can be constructed starting from the available a priori knowledge on the motor model (e.g. resistances and reactances are positive) and its choice is not critical for the identification procedure, as long as the chosen set is large enough to contain parameter values that achieve a good fitting between the model outputs and the identification data-set. As regards the constraints \eqref{E:opt_probl_eq}, they account for the model equations described in Section \ref{SS:model} and suitably discretized in order to numerically integrate them.

In this paper, we consider and compare different alternatives for the discrete-time model equations in \eqref{E:opt_probl_eq} and for the cost matrix $\boldsymbol{J}(\boldsymbol{p})$ in \eqref{E:opt_probl_cost}, as detailed in the following sub-sections. 

\subsection{Model discretization}
\label{SS:Discr_meth}
The induction motor model has to be discretized for the sake of numerical integration. Since the measurements coming from the sensor boxes and the smart circuit breaker are acquired at discrete time instants with sampling period $t_s$, we decide to integrate the model numerically with a fixed integration step equal to $t_s$. Albeit not strictly necessary (since one can in principle employ a smaller integration step and then consider, to compute the fitting errors, the model outputs at the same time instants when the experimental data have been sampled), this choice simplifies the identification procedure and its implementation on industrial hardware. In this work, we tested the performance and properties of the estimation algorithm at different sampling rates (and corresponding integration steps) and using two numerical integration techniques: the forward Euler method and an approach we called ''Input preview'' method.
\subsubsection{Forward Euler}
Discretizing the state equation \eqref{E:compactcontmatmodel} of the motor model with the forward Euler method yields:
\begin{equation}
\label{E:Eulermatmodel}
\hat{\boldsymbol{x}}(k+1)=\Big(\boldsymbol{I}+t_s\boldsymbol{A}(\hat{\omega}_r(k))\Big)\hat{\boldsymbol{x}}(k)+t_s\boldsymbol{B}\boldsymbol{u}(k)+t_s \boldsymbol{\beta}(\hat{\boldsymbol{x}}(k)),
\end{equation}
with $\hat{\boldsymbol{x}}(0)=0$. 
\subsubsection{Input preview}
The discrete-time expression of \eqref{E:compactcontmatmodel} obtained using the ''Input preview'' method is:
\begin{equation}
\label{E:Tustinmatmodel}
\hat{\boldsymbol{x}}(k+1)= \left(\boldsymbol{I}-\frac{t_s}{2}\boldsymbol{A}(\hat{\omega}_r(k))\right)^{-1} \left( \Big(\boldsymbol{I}+\frac{t_s}{2}\boldsymbol{A}(\hat{\omega}_r(k))\Big)\hat{\boldsymbol{x}}(k) + \frac{t_s}{2}\boldsymbol{B}\Big(\boldsymbol{u}(k+1)+\boldsymbol{u}(k)\Big)+t_s \boldsymbol{\beta}(\hat{\boldsymbol{x}}(k)) \right)
\end{equation}
with $\hat{\boldsymbol{x}}(0)=0$. The proposed discretization method is inspired in a sense by the Tustin approach, since it considers  the one-step-ahead input value $u(k+1)$, and the ``forward projection'' of the linear part of the system's dynamics, given by $\left(I- \frac{t_s}{2} \boldsymbol{A}(\omega_r(k+1)) \right)^{-1}$. The nonlinear part of matrix $\boldsymbol{A}$, which features the linear dependence of some coefficients on the rotor speed $\omega_r(t)$, can be approximated with a linear one by a timescale separation argument, since the rotor speed dynamics are significantly slower than the electric ones; thus, inside $\boldsymbol{A}$ we can consider $\omega_r(k+1) \simeq \omega_r(k)$, as done in \eqref{E:Tustinmatmodel}. We resort to this discretization method since, without an LTI system at hand, the Tustin approach does not allow us to obtain an explicit expression for $\hat{\boldsymbol{x}}(k+1)$. Therefore, we propose the Input preview approach which, as the experimental results presented in Section \ref{S:results} show, achieve better performance than the forward Euler method, while still retaining a reasonably low computational complexity (since it does not require an iterative numerical solution at each time step, like implicit integration methods do) and allowing one to derive an explicit calculation of parameter sensitivities, which we exploit to compute the cost function's gradient and estimate its Hessian. The latter aspect improves greatly  the computational efficiency when solving the identification problem. Finally, regarding the input vector $\boldsymbol{u}(k)$, this corresponds to the $k-$th column of either matrix $\tilde{\boldsymbol{V}}_{SB}$ (if the voltage measured by the sensor box is used) or $\tilde{\boldsymbol{V}}_{CB}$ (if the smart breaker data are used).
\subsubsection{Output equation}
Since they are static relationships, the discrete-time output equations are equal to the continuous-time ones. Thus, in case of stator currents (to be compared with the sensor boxes' measurements) we have:
\begin{equation}
\label{E:SB_out_disc}
\hat{\boldsymbol{y}}_{SB}(k)=\boldsymbol{C}\hat{\boldsymbol{x}}(k),
\end{equation}
where $\hat{\boldsymbol{x}}(k)$ is computed by iterating either \eqref{E:Eulermatmodel} or \eqref{E:Tustinmatmodel}. \\
\noindent In the case of stator current derivatives (to be compared with the circuit breaker  measurements) the discrete-time output equation is:
\begin{equation}
\label{E:CB_out_disc}
\hat{\boldsymbol{y}}_{CB}(k)=\boldsymbol{C}\dot{\hat{\boldsymbol{x}}}(k),
\end{equation}
where $\dot{\hat{\boldsymbol{x}}}(k)$ is defined as
\begin{equation}
\label{E:xdot_discr}
\dot{\hat{\boldsymbol{x}}}(k)=\boldsymbol{A}(\hat{\omega}_r(k))\hat{\boldsymbol{x}}(k)+\boldsymbol{B}\boldsymbol{u}(k)+\boldsymbol{\beta}(\hat{\boldsymbol{x}}(k)).
\end{equation}

Each one of the output equations \eqref{E:SB_out_disc} and \eqref{E:CB_out_disc} can be combined with either \eqref{E:Eulermatmodel} or \eqref{E:Tustinmatmodel} and inserted in the constraints \eqref{E:opt_probl_eq} to obtain four possible cases, i.e. Euler or Input preview discretization and either currents or current derivatives as model outputs. In the literature, to the best of our knowledge only the case of Euler integration and current outputs has been considered so far, while here we explore all four combinations. Moreover, each case can be considered with a different available sampling frequency. In our experimental results we provide a thorough analysis of the results obtained with all these alternatives.

\subsection{Cost function definition}
The cost matrix $\boldsymbol{J}(\boldsymbol{p})$ in \eqref{E:opt_probl_cost} changes depending on whether current or current derivative data are employed for the fitting criterion. In the following, we assume that a data-set to be used for the parameter identification phase has been fixed, such that the cost matrix depends only on the parameter value $\boldsymbol{p}$. 
\subsubsection{Current data}
In the case of current data (i.e. acquired by the sensor boxes in our setup), the cost matrix is computed as
\begin{equation}
\label{E:J_sb}
\boldsymbol{J}_{SB}(\boldsymbol{p})= \left( \left( \tilde{\boldsymbol{I}}-\hat{\boldsymbol{Y}}(\tilde{\boldsymbol{V}}_{SB},\boldsymbol{p}) \right) \left( \tilde{\boldsymbol{I}}-\hat{\boldsymbol{Y}}(\tilde{\boldsymbol{V}}_{SB},\boldsymbol{p}) \right)^T \right)
\end{equation}
where $\hat{\boldsymbol{Y}}(\tilde{\boldsymbol{V}}_{SB},\boldsymbol{p}):=[\hat{\boldsymbol{y}}_{SB}(1),\cdots,\hat{\boldsymbol{y}}_{SB}(N)]\in\mathbb{R}^{2\times N}$ 
is a matrix containing the stator current signals in the \textit{dq}-components simulated with the motor model (either \eqref{E:Eulermatmodel} or \eqref{E:Tustinmatmodel}) and the output equation \eqref{E:SB_out_disc}, excited by the stator voltage signal $\tilde{\boldsymbol{V}}_{SB}$ as input.
\subsubsection{Current derivative data}
In the case of current derivative measurements (i.e. acquired by the smart breaker), the objective function is defined as
\begin{equation}
\label{E:J_emax}
\boldsymbol{J}_{CB}(\boldsymbol{p})=\left( \left( \tilde{\dot{\boldsymbol{I}}}-\hat{\dot{\boldsymbol{Y}}}(\tilde{\boldsymbol{V}}_{CB},\boldsymbol{p}) \right) \left( \tilde{\dot{\boldsymbol{I}}}-\hat{\dot{\boldsymbol{Y}}}(\tilde{\boldsymbol{V}}_{CB},\boldsymbol{p}) \right)^T \right)
\end{equation}
where $\hat{\dot{\boldsymbol{Y}}}(\tilde{\boldsymbol{V}}_{CB},\boldsymbol{p}):=[\hat{\boldsymbol{y}}_{CB}(1),\cdots,\hat{\boldsymbol{y}}_{CB}(N)]$ contains the simulated stator current derivative signals in the \textit{dq}-components, obtained by integrating the motor model  with  the output equation \eqref{E:CB_out_disc} and applying the measured voltage sequence $\tilde{\boldsymbol{V}}_{CB}$ as input. 


\subsection{Algorithmic aspects: explicit gradient computation and Hessian estimates, parameter initialization}\label{SS:algo}
We solve the optimization problem \eqref{E:opt_probl} with a constrained Gauss-Newton algorithm, which is a sequential quadratic programming (SQP) approach where the computation of the gradient and the estimation of the Hessian of the cost function (required for the quadratic sub-problems) are computed by exploiting the problem structure \cite{NoWr06}. In particular, $J(\boldsymbol{p})$ can be re-written as
\begin{equation}\label{E:cost_F}
\boldsymbol{J}(\boldsymbol{p})=\boldsymbol{F}(\boldsymbol{p})^T \boldsymbol{F}(\boldsymbol{p}),
\end{equation} 
where $\boldsymbol{F}(\boldsymbol{p})\in\mathbb{R}^{2N}$ is a vector containing the differences between the measured outputs (currents or their derivatives) and the model outputs at each time step. The Jacobian matrix $\nabla_{\boldsymbol{F}}(\boldsymbol{p})$ of $\boldsymbol{F}(\boldsymbol{p})$ can be obtained by differentiating the model equations, resulting in a recursive formulation that can be computed together with the model integration. Such a recursion  is detailed in the Appendix for both Euler and Input preview discretization methods. Then, the gradient of $\boldsymbol{J}(\boldsymbol{p})$ is computed as:
\[
\nabla_{\boldsymbol{J}}(\boldsymbol{p})=2\nabla_{\boldsymbol{F}}(\boldsymbol{p})^T \boldsymbol{F}(\boldsymbol{p}),
\]
and, by considering the Taylor expansion of $\boldsymbol{F}(\boldsymbol{p})$ truncated at the first order and inserting it in \eqref{E:cost_F}, the Hessian of $\boldsymbol{J}(\boldsymbol{p})$  is approximated as 
\[
\nabla^2_{\boldsymbol{J}}(\boldsymbol{p})\simeq \nabla_{\boldsymbol{F}}(\boldsymbol{p})^T \nabla_{\boldsymbol{F}}(\boldsymbol{p}).
\]
In our experiments, the use of this approach resulted in significant computational savings with respect to general-purpose nonlinear programming solvers, as we mention in Section \eqref{S:results}.

Another relevant aspect from the point of view of computational efficiency is the initialization of the optimization routine. Due to the non-convex nature of the problem, the algorithm generally converges to a local optimum. The sub-optimality of the solution and the convergence speed are clearly sensitive to the initialization of the parameters required by the SQP solver.  One possible heuristic approach to attain  the global optimum is to run several times the algorithm with different initialization values,  generated randomly within the set of possible parameters $\mathcal{P}$, and then to consider the estimate that provides the minimum prediction error calculated over the batch of data. We adopt this approach here, and in Section \ref{SS:P_init} we analyze how the algorithm sensitivity to initialization changes with different sampling frequencies and discretization methods.
\\
\section{Experimental results}
\label{S:results}
We applied the described estimation procedure to real data acquired with the experimental rig described in Section \ref{S:probl_setup_descr}, consisting of about 100 direct on-line start-up transients of the 3-phase induction motors in different combinations (e.g. only motor $M1$, only $M2$, and various combined start-ups with different delays between the two motors). 

We considered different sampling frequencies and, for each one, we employed the data from one start-up experiment for the identification, and the data of three additional experiments for validation. In particular, as performance metric to compare the different tests we employ the Normalized Mean Prediction Error (NMPE), defined as 
\begin{equation*}
\textnormal{NMPE}:= \sqrt{\frac{tr\left(\left( \tilde{\boldsymbol{I}}-\hat{\boldsymbol{Y}}(\tilde{\boldsymbol{V}}_{SB},\hat{\boldsymbol{p}}) \right)\left( \tilde{\boldsymbol{I}}-\hat{\boldsymbol{Y}}(\tilde{\boldsymbol{V}}_{SB},\hat{\boldsymbol{p}}) \right)^T\right)}{tr \left( \tilde{\boldsymbol{I}} \tilde{\boldsymbol{I}}^T \right)}}.
\end{equation*} 
Note that in the NMPE calculation we always consider the motor currents as measured by the high-quality sensors installed in the sensor boxes. This means that also the parameters estimated from the smart breaker data (i.e. using current derivatives as identification data-set) are then tested against current measures collected by the sensor boxes. In all the results presented in the following, we provide the range of NMPE given by the minimum and maximum values obtained in the three validation experiments related to each specific test case.

As regards the set of admissible parameters $\mathcal{P}$, we selected rather wide ranges for each of the parameters to be estimated:
$$
\begin{bmatrix}
0 \\ 0 \\ 0 \\ 0 \\ 0 \\ 0 \\ 0
\end{bmatrix} \preceq \boldsymbol{p} \preceq \begin{bmatrix}
100 \\ 100 \\ 100 \\ 500 \\ 20 \\ 100 \\ 0.35
\end{bmatrix}
$$

In all the tests reported in the following, we employed the SQP solver based on the constrained Gauss-Newton approach and analytic computation of the gradient, as described in Section \ref{SS:algo}. The solver, implemented in MatLab, was able to converge on average in about 20 iterations and 120 s on a Laptop equipped with Intel i7 dual-core processor with 2.4 GHz clock speed and 8 GB of RAM. For a comparison, on the same hardware a standard optimization routine (MatLab \verb|fmincon|) took on average 50 iterations and 2200 s with the same termination tolerances.

\subsection{Sensors comparison}
\label{SS:senscomp}
To determine whether the data collected by the industrial voltage and current sensors installed in the considered commercial circuit breaker are good enough to identify the model parameters, we compared the results of the estimation procedure performed using data acquired by the sensor boxes, which have a maximum sampling frequency 5 kHz, with the results obtained using data measured by the smart breaker, where we selected a sampling frequency of 4.8 kHz. In both cases, the discrete-time model is obtained using the Input preview method. The results related to motor  $M1$ are presented in Table \ref{T:SboxvsEmax}. It can be noted that the differences between the two parameter estimates and the resulting NMPE ranges are minimal. 

\begin{table}[thpb]
\caption{Motor $M1$ - Identified parameters: comparison between the results obtained with sensor box data and circuit breaker data.}
\label{T:SboxvsEmax}
\centering
\begin{adjustbox}{max width=1\columnwidth}
\begin{tabular}{| c | c | c |}
%
\hline
 & Sensor box & Circuit breaker \\ \hline
$R_s$ & 0.48 & 0.48 \\ \hline
$R_r$ & 0.20 & 0.21 \\ \hline
$X_l$ & 0.29 & 0.30 \\ \hline
$X_m$ & 11.92 & 11.29 \\ \hline
$J_r$ & 0.26 & 0.26 \\ \hline
$T_{l_0}$ & 0 & 0 \\ \hline
$T_{l_1}$ & 0.039 & 0.037 \\ \hline 
NMPE & 0.0810 $\div$ 0.0813 & 0.0775 $\div$ 0.0798 \\ \hline
\end{tabular}
\end{adjustbox}
\end{table}
\begin{figure*}[thpb]
      \centering
      		\begin{tabular}{cc}
      			 (a) & (b) \\
     			 \includegraphics[width=0.48\columnwidth]{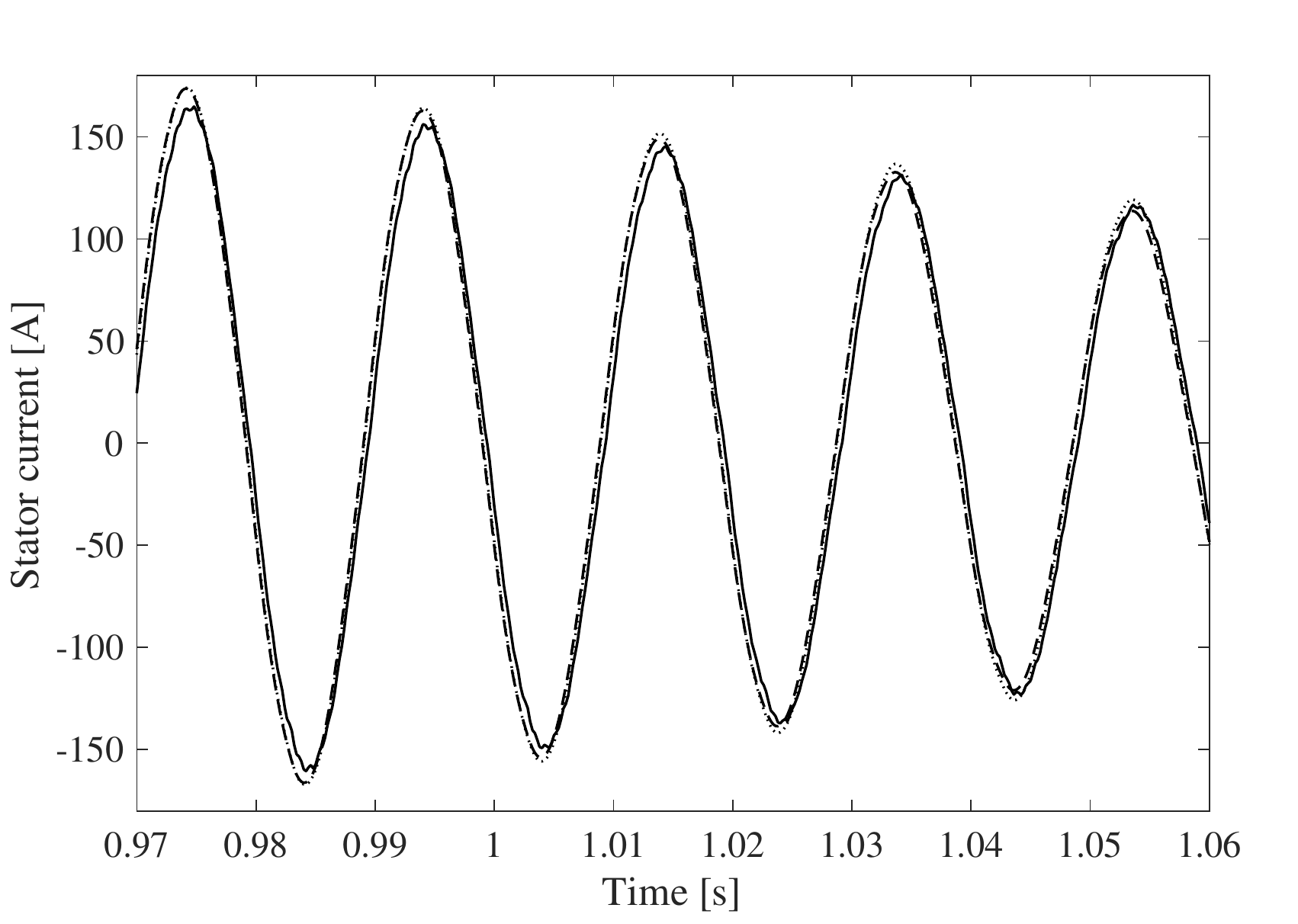} & \includegraphics[width=0.48\columnwidth]{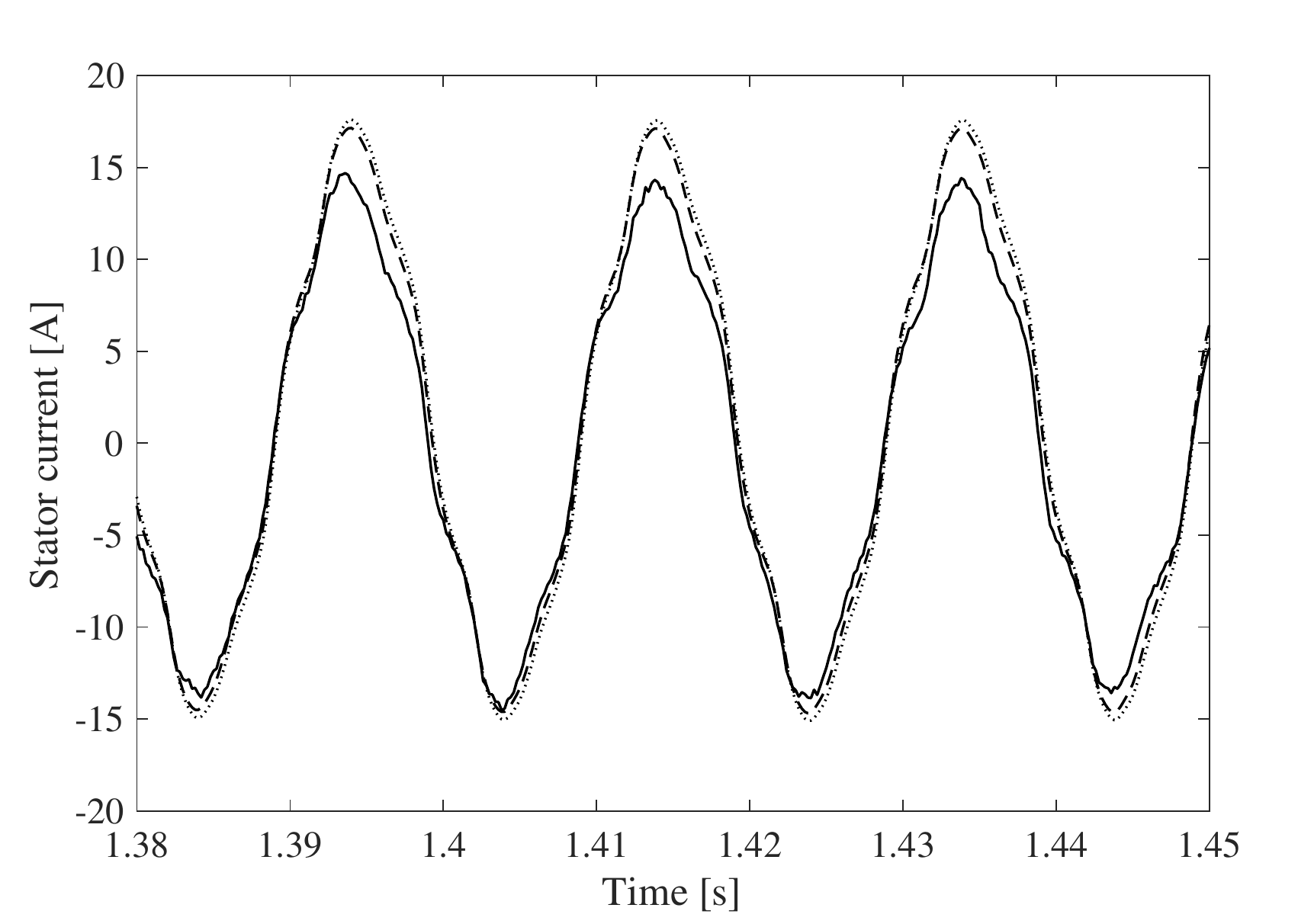} \\ \\
     			 (c) & (d) \\
     			 \includegraphics[width=0.48\columnwidth]{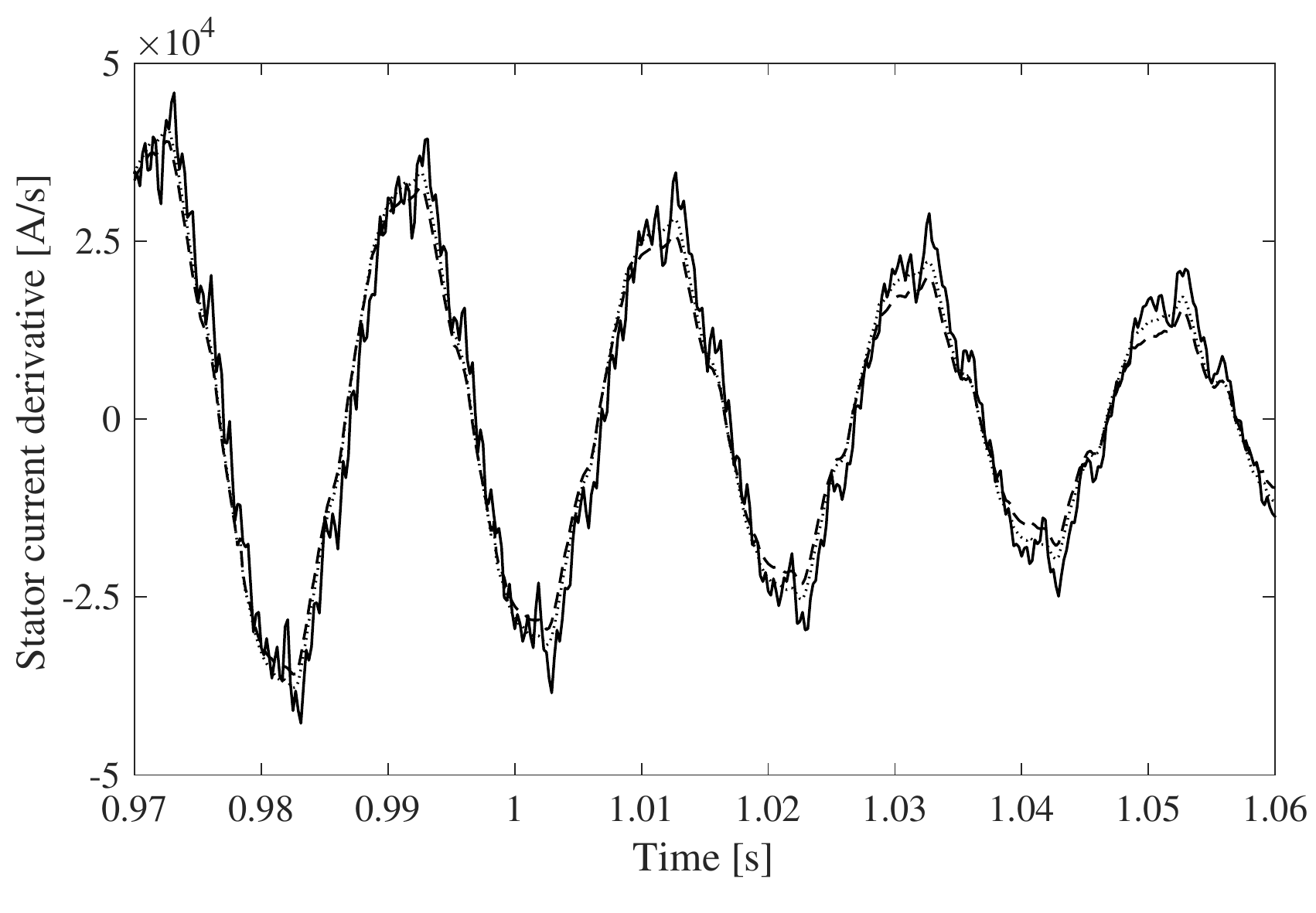} & \includegraphics[width=0.48\columnwidth]{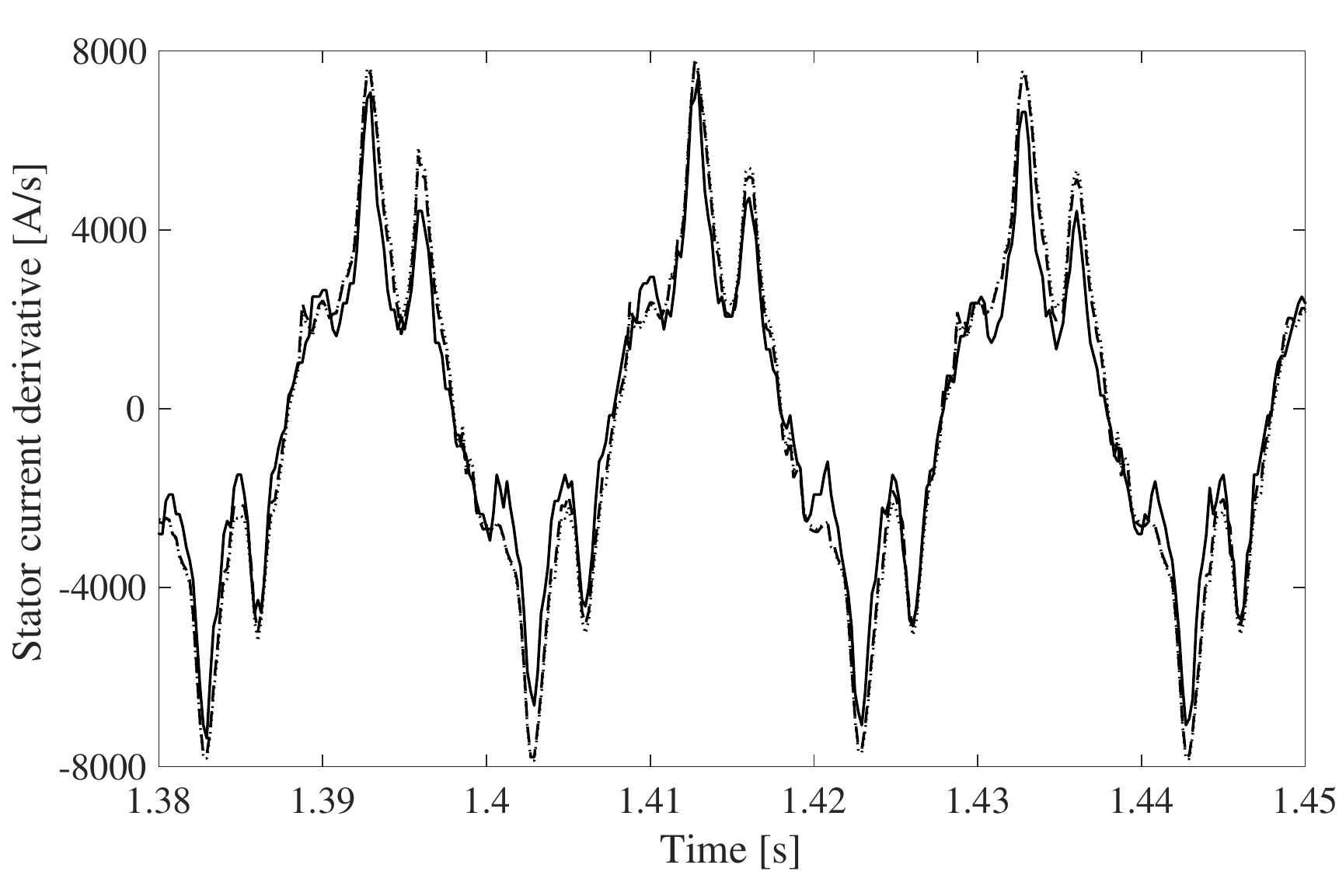}
      		\end{tabular}
      \caption{Motor $M1$ - Estimated parameters validation. Comparison between simulated and measured current signals acquired by sensor box (plots (a) and (b)), and between simulated and measured current derivative signals acquired by the smart circuit breaker (plots (c) and (d)). The plots on the left pertain to the first part of the transient, while those on the right pertain to a steady speed condition. Solid lines: measured $q$ component of stator current (or current derivative); dashed lines: simulated $q$ component based on parameters estimated from sensor box data: dotted lines: simulated $q$ component based on parameters estimated from circuit breaker data.}
      \label{f:caso1_SBvsEMAX}
\end{figure*} 
Figs. \ref{f:caso1_SBvsEMAX}(a) and \ref{f:caso1_SBvsEMAX}(b) show the comparison between the $q$ component of the stator current, measured by the sensor box during a validation experiment, and the signal reconstructed using the parameters identified from the sensor box data-set and from the circuit breaker data-set. The fitting is good in both cases, as expected from the NMPE results of Table \ref{T:SboxvsEmax}.
Figs. \ref{f:caso1_SBvsEMAX}(c) and \ref{f:caso1_SBvsEMAX}(d) present the comparison between the $q$ component of the stator current derivatives, measured by the smart circuit breaker during the validation experiment, and the current derivatives reconstructed using the parameters identified from the two different data-sets described before. Also in this case the fitting is good and the parameters estimated from the two different data-sets have a comparable performance. \\
Table \ref{T:SboxvsEmax_M2} presents the comparison between the parameters of motor $M2$ estimated from sensor box data and from smart circuit breaker data. The obtained results are fully aligned with those of motor $M1$.

These results indicate that it is possible to obtain a good estimate of the motor parameters also using data acquired by Rogowski coil sensors from commercial circuit breakers, where the stator current derivatives are measured in place of a direct measure of the stator currents. In the remainder of this section, we investigate more in detail the performance of the estimation algorithm using the data acquired by the smart circuit breaker, with different choices of discretization method and sampling frequency, and we analyze the sensitivity to parameter initialization and model over-parametrization.
\begin{table}[thpb]
\caption{Motor $M2$ - Identified parameters: comparison between the results obtained with sensor box data and circuit breaker data.}
\label{T:SboxvsEmax_M2}
\centering
\begin{adjustbox}{max width=1\columnwidth}
\begin{tabular}{| c | c | c |}
%
\hline
 & Sensor box & Circuit breaker \\ \hline
$R_s$ & 1.15 & 1.12 \\ \hline
$R_r$ & 0.49 & 0.49 \\ \hline
$X_l$ & 0.70 & 0.71 \\ \hline
$X_m$ & 33.67 & 34.48 \\ \hline
$J_r$ & 0.27 & 0.28 \\ \hline
$T_{l_0}$ & 0 & 0 \\ \hline
$T_{l_1}$ & 0.035 & 0.031 \\ \hline 
NMPE & 0.0853 $\div$ 0.0882 & 0.0954 $\div$ 0.0976 \\ \hline
\end{tabular}
\end{adjustbox}
\end{table}
\subsection{Comparison between discretization methods}
\label{SS:discrmeth}
We applied the estimation algorithms derived using the two discretization methods described in Section \ref{SS:Discr_meth} to data-sets acquired by the circuit breaker with various sampling frequencies. 
Table \ref{T:TustvsEA} contains the parameter values identified using  one data-set acquired at $f_s=4.8$ kHz, and another one at $f_s=2.4$ kHz. In the Table, we highlight in bold the parameter values that are clearly different from the best ones, reported previously in Section \ref{SS:senscomp}.  

For the case of the data-set acquired with $f_s=4.8$ kHz, the results of the estimation procedure based on the Euler and the Input preview methods are similar, and the estimation errors resulting from the two cases are comparable. On the other hand, the data-set acquired with $f_s=2.4$ kHz leads to significantly different results: the estimates obtained with the forward Euler method are not consistent with those obtained by the same method at higher frequency, and the NMPE values are much larger. The estimates obtained with the Input preview method appear to be resilient to lower frequencies and they are still very close to the best ones. This is coherent with the system discretization theory, as the trapezoidal integration method gives more precise and stable results than the forward Euler method.

\begin{figure}[thpb]
      \centering
      \includegraphics[width=0.8\columnwidth]{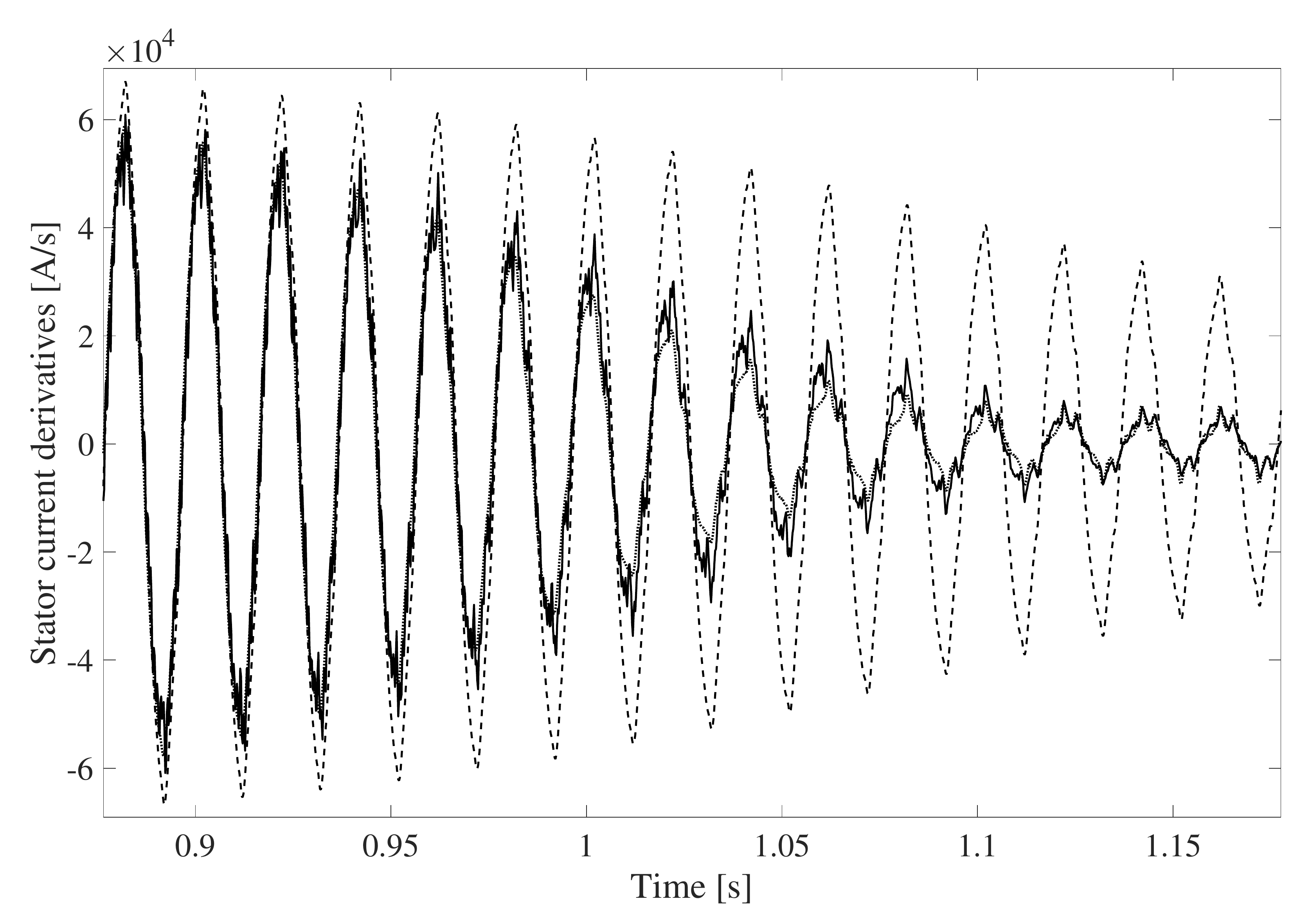}
      \caption{Motor $M1$ - validation of parameter estimates obtained with circuit breaker data measured at $f_s=2.4$ kHz and different discretization methods. The course of the $q$ component of the stator current derivative is shown. Solid line: measured signal at $f_s=4.8$ kHz; dashed line: simulated signal obtained with parameters identified using the forward Euler method; dotted line: simulated signal obtained with parameters identified using the Input preview method.}
      \label{f:EA-Tust_freq}
\end{figure}

We also cross-validated the results of the estimation procedure, by comparing the current derivative signals acquired at $f_s=4.8$ kHz with the ones simulated using the parameters estimated from a data-set measured at a different sampling frequency, e.g. $f_s=2.4$ kHz. An example of the results is reported in Fig. \ref{f:EA-Tust_freq}: it is clear that the algorithm based on the forward Euler method gives worse results compared to the Input preview case. Looking at the simulated current derivative reported in Fig. \ref{f:EA-Tust_freq}, it is clearly visible that in the forward Euler case both the transient and the steady-state part of the motor start-up do not match with the measured signal, while in the Input preview case the fitting is still fairly good. 
\begin{table}[thpb]
\caption{Motor $M1$ - Identified parameters: comparison between different discretization methods. We highlight with bold characters the parameters that are considerably different from the optimum ones.}
\label{T:TustvsEA}
\centering
\begin{adjustbox}{max width=1\columnwidth}
\begin{tabular}{| c | c | c |}
\hline
 & \multicolumn{2}{c |}{Input preview} \\ \hline
 & 4.8 kHz & 2.4 kHz \\ \hline
$R_s$ & 0.48 & 0.48 \\ \hline
$R_r$ & 0.21 & 0.21 \\ \hline
$X_l$ & 0.30 & 0.30 \\ \hline
$X_m$ & 11.29 & 11.28 \\ \hline
$J_r$ & 0.26 & 0.26 \\ \hline
$T_{l_0}$ & 0 & 0 \\ \hline
$T_{l_1}$ & 0.037 & 0.037 \\ \hline 
NMPE & 0.0775 $\div$ 0.0798 & 0.0784 $\div$ 0.0791 \\ \hline 
\end{tabular} \quad \begin{tabular}{| c | c | c |}
\hline
 & \multicolumn{2}{c |}{Euler} \\ \hline
 & 4.8 kHz & 2.4 kHz \\ \hline
$R_s$ & 0.53 & 0.55 \\ \hline
$R_r$ & 0.22 & 0.24 \\ \hline
$X_l$ & 0.28 & 0.26 \\ \hline
$X_m$ & \textbf{2.52} & \textbf{1.45} \\ \hline
$J_r$ & 0.24 & 0.19 \\ \hline
$T_{l_0}$ & 0 & \textbf{13.21} \\ \hline
$T_{l_1}$& \textbf{0.13} & \textbf{0.20} \\ \hline 
NMPE & 0.1148 $\div$ 0.1162 & \textbf{0.3268} $\div$ \textbf{0.3337} \\ \hline
\end{tabular}
\end{adjustbox}
\end{table}

\begin{figure}[thpb]
      \centering
      \includegraphics[width=0.8\columnwidth]{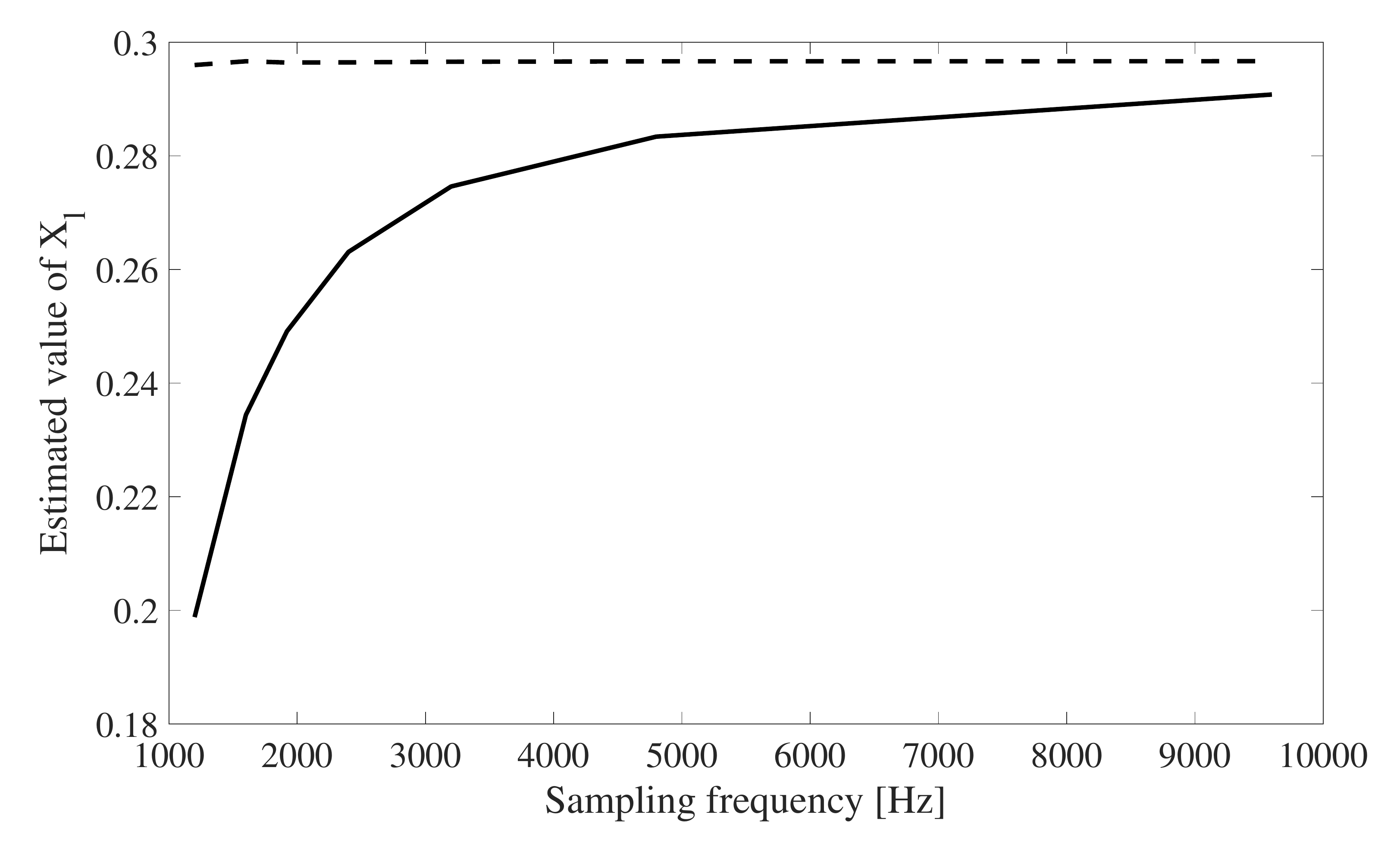}
      \caption{Motor $M1$ - circuit breaker data, sensitivity to data sampling frequency. Identified value of the reactance $X_l$ as a function of the data sampling frequency: estimation algorithm based on the forward Euler method (solid line) and on the Input preview method (dashed).}
      \label{f:Xl_EA-Tust_freq}
\end{figure}

\begin{figure}[thpb]
      \centering
      \includegraphics[width=0.8\columnwidth]{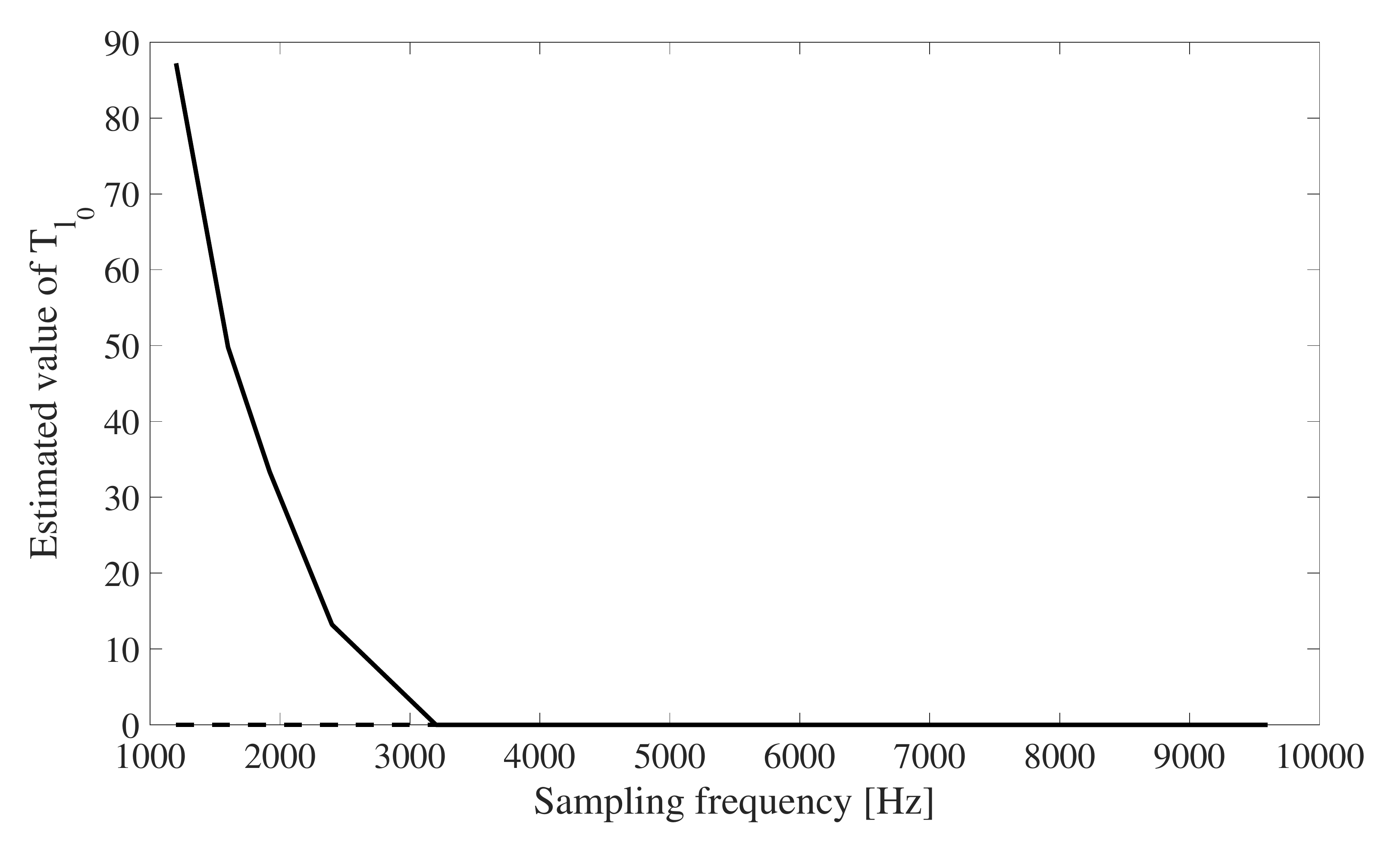}
      \caption{Motor $M1$ - circuit breaker data, sensitivity to data sampling frequency. Identified value of the constant load torque parameter $T_{l_0}$ as a function of the data sampling frequency: estimation algorithm based on the forward Euler method (solid line)and on the Input preview method (dashed).}
      \label{f:Tl0_EA-Tust_freq}
\end{figure}

We performed a more detailed analysis of the sensitivity of the estimation results to the sampling frequency, by running the estimation algorithm on data acquired with $f_s$ spacing from $1.2$ kHz to $9.6$ kHz. From these experiments, whose results are partially depicted in Figs. \ref{f:Xl_EA-Tust_freq} and \ref{f:Tl0_EA-Tust_freq}, it is possible to see that, in the case of  forward Euler method, the identified parameters vary sensibly with the sampling frequency, while, in the Input preview case, they exhibit a much lower variability.

The reported results further confirm that the Input preview method is generally more stable with respect to variations of the sampling frequency of the measured data (and size of the integration step), while the forward Euler estimate diverges at low frequencies.
\subsection{Sensitivity to parameter initialization}
\label{SS:P_init}
\begin{figure*}[thpb]
      \centering
      		\begin{tabular}{cc}
      			 (a) & (b) \\
     			 \includegraphics[width=0.48\columnwidth]{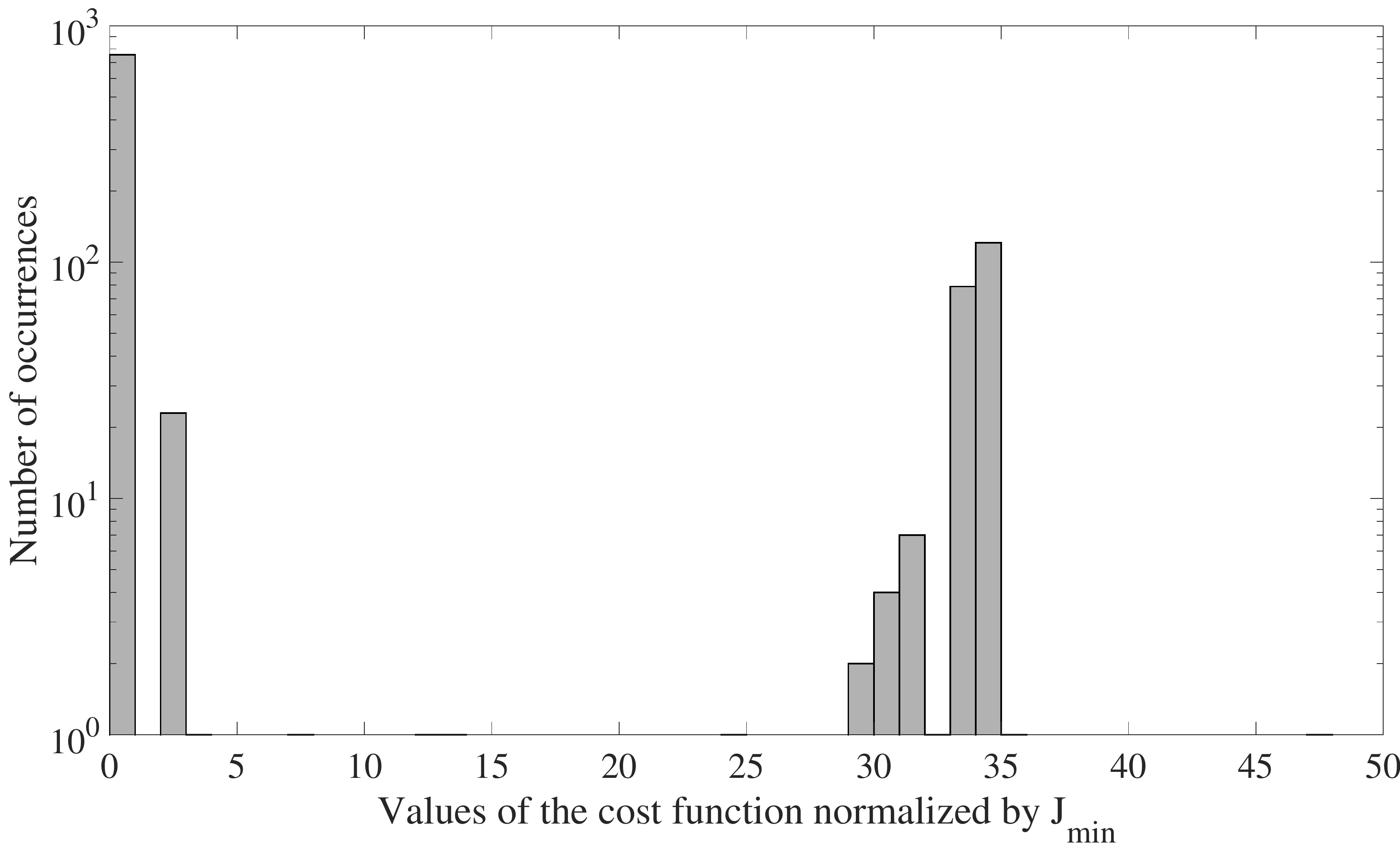} & \includegraphics[width=0.48\columnwidth]{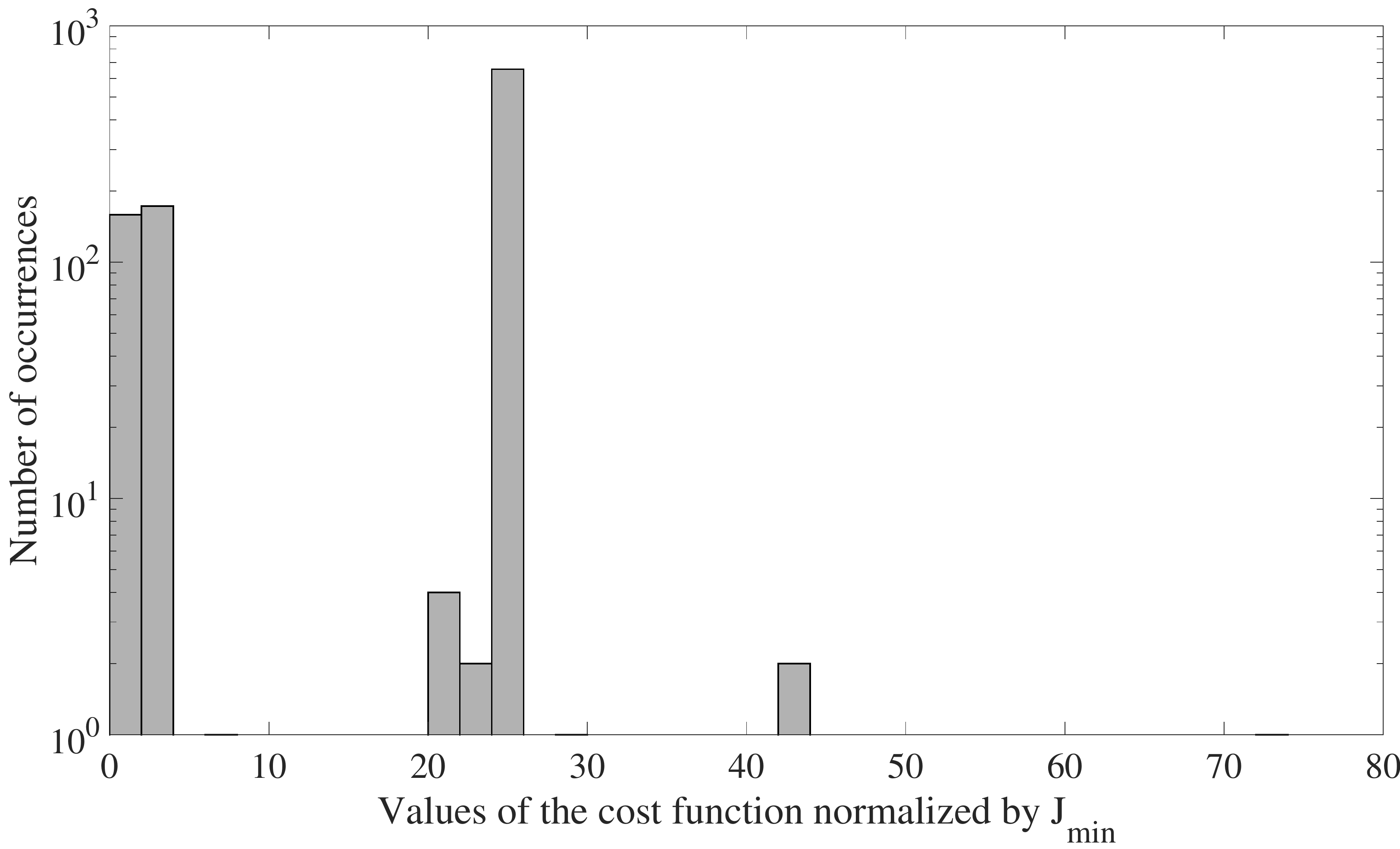} \\ \\
     			 (c) & (d) \\
     			 \includegraphics[width=0.48\columnwidth]{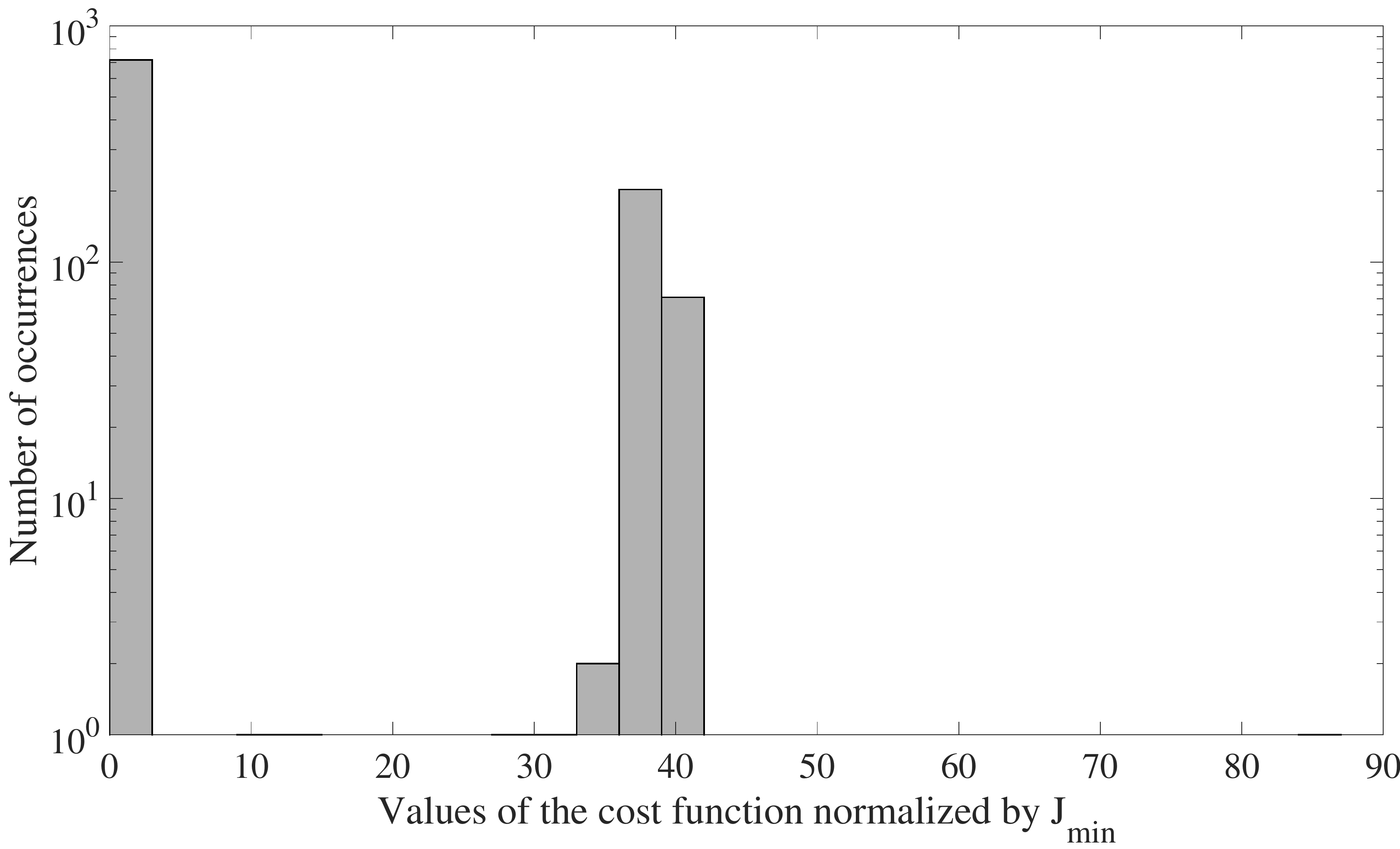} & \includegraphics[width=0.48\columnwidth]{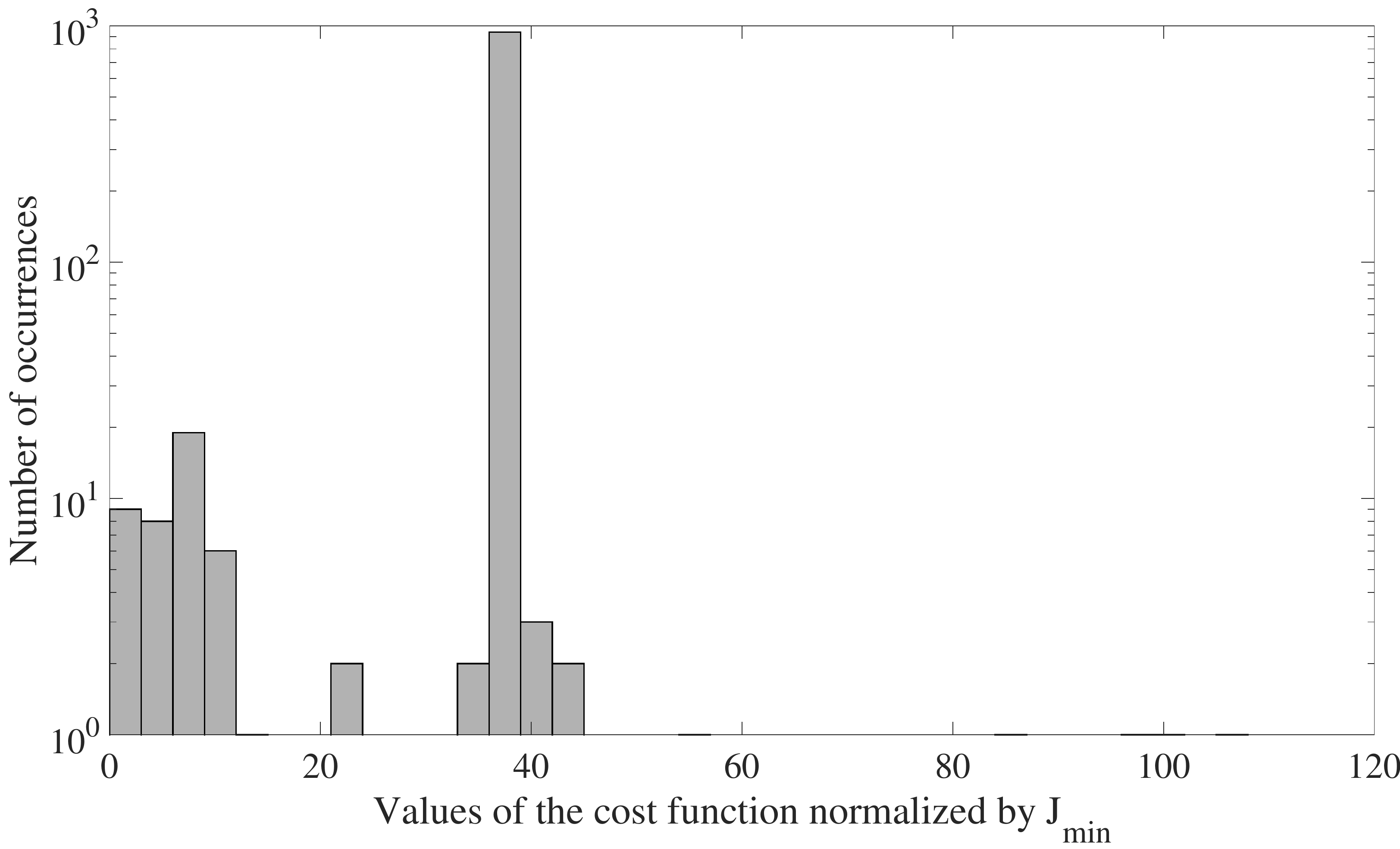}
      		\end{tabular}
      \caption{Motor $M1$ - Convergence analysis of 1000 estimation routines with uniformly distributed random-generated initial parameters. Identification data measured by the circuit breaker. Discretization method: Input preview: (a) data acquired at 4.8 kHz; (b)  data acquired at 2.4 kHz. Discretization method: forward Euler: (c) data acquired at 4.8 kHz; (d)  data acquired at 2.4 kHz.}
      \label{f:Istogrammi_T_e_EA}
\end{figure*}

As stated in Section \ref{S:nonlinid}, due to the non-convex nature of the estimation problem the optimality of the computed estimates depends on the initialization of the optimization algorithm. To study the sensitivity of the algorithm to initialization, we performed 1000 estimation routines on the same data-set, where we randomly picked $\boldsymbol{p}_0$, the initial parameter value, with a uniform distribution over the following sub-set of $\mathcal{P}$:
$$
\begin{bmatrix}
0 \\ 0 \\ 0 \\ 0 \\ 0 \\ 0 \\ 0
\end{bmatrix} \preceq \boldsymbol{p}_0 \preceq \begin{bmatrix}
10 \\ 10 \\ 10 \\ 15 \\ 2 \\ 1 \\ 0.042
\end{bmatrix}.$$
Note that the value of $\boldsymbol{p}$ that we consider as the global optimum is inside the described set (compare Table \ref{T:SboxvsEmax}).\\
We repeated this analysis using data-sets acquired at different frequencies, in order to understand how the sensitivity to the parameters initialization varies with respect to $f_s$, and using both discretization methods. In particular, we ran the 1000 estimation routines for each of the possible combinations of sampling frequency (2.4 kHz or 4.8 kHz) and discretization method. We considered as acceptable the result of an estimation routine only when the value of the cost function $\boldsymbol{J}(\boldsymbol{p})$, evaluated at the obtained $\hat{\boldsymbol{p}}$, is inside the interval $\left[ \boldsymbol{J}_{min}, \; 1.05\cdot \boldsymbol{J}_{min} \right]$, where $\boldsymbol{J}_{min}$ is the minimum value of $\boldsymbol{J}(\boldsymbol{p})$ obtained across all 1000 tests for the specific combination of sampling frequency and discretization method. The results of this analysis are reported in Table \ref{T:randP0}. 
\begin{table}[thpb]
\caption{Motor $M1$ - Number of acceptable estimation results for uniformly distributed random-generated initial parameters.}
\label{T:randP0}
\centering
\begin{tabular}{| c | c | c |}
\hline
 & Input preview & Euler \\ \hline
$f_s=4.8$ kHz & 756 over 1000 & 703 over 1000 \\ \hline
$f_s=2.4$ kHz & 159 over 1000 & 7 over 1000 \\ \hline
\end{tabular}
\end{table}
\\
\noindent Fig. \ref{f:Istogrammi_T_e_EA} illustrates the distribution of the final values of the objective function (i.e. at the termination of each identification procedure) obtained by the estimation routines, normalizing each result by the corresponding value of $\boldsymbol{J}_{min}$ considered as the optimum. These figures show how many times the estimation algorithm reaches the best achieved fitting cost and how the results distribute around  local minima with different degree of sub-optimality. 
%

The obtained results indicate a clear increase of the sensitivity to the initialization as the sampling frequency decreases, for both discretization methods. Moreover, the estimation algorithm based on the forward Euler method exhibits a greater sensitivity than the one based on trapezoid approximation.
\subsection{Effects of over-parametrization}
\label{SS:overparam}
In this work we chose to adopt a general purpose model of the load torque of the induction motor, that is based on two parameters; in the experimental setup, as described in Section \ref{S:probl_setup_descr}, the load torque applied to the motor is only proportional to its angular speed, without any constant term. We decided to analyze how the estimation results vary if we use a model of the motor load that is based only on the linear viscous term, which corresponds to the physical behavior of our experimental setup. In Table \ref{T:Tl0siono} we compare the results of the estimation algorithm applied to data acquired with different sampling frequencies, in the case of the over-parametrized load torque model and in the case having only the linear viscous term. From this results, it is clear that the Input preview method is sufficiently robust to provide good performance also with a generic and over-parametrized load model; on the other side, we can confirm that the forward Euler method gives less consistent results as the sampling frequency lowers, with both the load torque models.
\begin{table}[thpb]
\caption{Motor M1 - Identified parameters: comparison between different load models. We highlight with bold characters the parameters that are considerably different from the optimum ones.}
\label{T:Tl0siono}
\centering
\begin{adjustbox}{max width=1\columnwidth}
\begin{tabular}{| c | c | c | c | c |}
\hline
 & \multicolumn{4}{c |}{Input preview} \\ \hline
 & 4.8 kHz & 4.8 kHz & 2.4 kHz & 2.4 kHz \\ 
 & with $T_{l_0}$ & without $T_{l_0}$ & with $T_{l_0}$ & without $T_{l_0}$ \\ \hline
$R_s$ & 0.48 & 0.48 & 0.48 & 0.48 \\ \hline
$R_r$ & 0.21 & 0.21 & 0.21 & 0.21 \\ \hline
$X_l$ & 0.30 & 0.30 & 0.30 & 0.30 \\ \hline 
$X_m$ & 11.29 & 11.29 & 11.28 & 11.28 \\ \hline 
$J_r$ & 0.26 & 0.26 & 0.26 & 0.26 \\ \hline 
$T_{l_0}$ & 0 & - & 0 & - \\ \hline 
$T_{l_1}$ & 0.037 & 0.037 & 0.037 & 0.037 \\ \hline 
\end{tabular} \quad \begin{tabular}{| c | c | c | c | c |}
\hline
 & \multicolumn{4}{c |}{Euler} \\ \hline
 & 4.8 kHz & 4.8 kHz & 2.4 kHz & 2.4 kHz \\ 
 & with $T_{l_0}$ & without $T_{l_0}$ & with $T_{l_0}$ & without $T_{l_0}$ \\ \hline
$R_s$ & 0.53 & 0.53 & 0.55 & 0.57 \\ \hline
$R_r$ & 0.22 & 0.22 & 0.24 & 0.22 \\ \hline
$X_l$ & 0.28 & 0.28 & 0.26 & 0.26 \\ \hline 
$X_m$ & \textbf{2.52} & \textbf{2.52} & \textbf{1.45} & \textbf{1.36} \\ \hline 
$J_r$ & 0.24 & 0.24 & 0.19 & 0.2 \\ \hline 
$T_{l_0}$ & 0 & - & \textbf{13.21} & - \\ \hline 
$T_{l_1}$ & \textbf{0.13} & \textbf{0.13} & \textbf{0.20} & \textbf{0.24} \\ \hline
\end{tabular}
\end{adjustbox}
\end{table}
\section{Conclusion and future directions}
\label{S:conclusions}
We presented an experimental study on the use of smart circuit breaker to identify the parameters of three-phase induction motors, using the data collected during direct-on-line start-ups. The main finding is that commercial circuit breakers are able to collect data with suitable quality to effectively carry out the identification procedure, as compared with high-performance laboratory sensors. We further investigated different variants of the identification algorithm and highlighted their advantages and drawbacks with substantial experimental data. Future research will be aimed to employ the identified models to carry out motor detection and monitoring tasks, as well as to achieve advanced protection performance by better discriminating between faults and motor inrush currents.
\appendix
\subsection*{Appendix}
\subsubsection*{Recursive equations for the computation of the cost function gradient.} 
We start with the recursive equations used to calculate the cost function gradient and Hessian in the forward Euler case. For a given value of $\boldsymbol{p}$, in the case of sensor box data (i.e. current output) we have $\boldsymbol{F}=\left[F_{SB}(1)^T,\cdots,F_{SB}(N)^T\right]^T$, with $F_{SB}(k)=\tilde{\boldsymbol{I}}_{SB}(k)-\hat{\boldsymbol{I}}_{SB}(k)$. The Jacobian $\nabla_{\boldsymbol{F}}(\boldsymbol{p})$ can be then obtained by differentiating $F_{SB}(k)=\tilde{\boldsymbol{I}}_{SB}(k)-\hat{\boldsymbol{I}}_{SB}(k)$ with respect to $\boldsymbol{p}$, resulting in the following recursive equations:
\begin{gather*}
\frac{\partial F_{SB}(k)}{\partial \boldsymbol{p}}=\frac{\partial \hat{\boldsymbol{I}}_{SB}(k)}{\partial \boldsymbol{p}}= \boldsymbol{C} \frac{\partial \hat{\boldsymbol{x}}(k)}{\partial \boldsymbol{p}}+\frac{\partial (\boldsymbol{C}\boldsymbol{x})}{\partial \boldsymbol{p}} \bigg\rvert_{\substack{\boldsymbol{x}=\hat{\boldsymbol{x}}(k)}}, \\
 \\
\frac{\partial \hat{\boldsymbol{x}}(k+1)}{\partial \boldsymbol{p}}= \frac{\partial \Big( \big(\boldsymbol{I}+t_s\boldsymbol{A}(\boldsymbol{x})\big)\boldsymbol{x} +t_s \boldsymbol{\beta}(\boldsymbol{x}) \Big)}{\partial \boldsymbol{x}}\bigg\rvert_{\substack{\boldsymbol{x}=\hat{\boldsymbol{x}}(k)}} \cdot \frac{\partial \hat{\boldsymbol{x}}(k)}{\partial \boldsymbol{p}} + \frac{\partial \Big( \big(\boldsymbol{I}+t_s\boldsymbol{A}(\boldsymbol{x})\big)\boldsymbol{x} +t_s \boldsymbol{\beta}(\boldsymbol{x}) \Big)}{\partial \boldsymbol{p}}\bigg\rvert_{\substack{\boldsymbol{x}=\hat{\boldsymbol{x}}(k)}}.
\end{gather*}  
In the case of smart circuit breaker data (i.e. current derivatives as output), by differentiating $F_{CB}(t)=\tilde{\dot{\boldsymbol{I}}}_{CB}(k)-\hat{\dot{\boldsymbol{I}}}_{CB}(k)$, we can obtain the following recursive equations:
\begin{gather*}
\frac{\partial F_{CB}(k)}{\partial \boldsymbol{p}}=\frac{\partial \hat{\dot{\boldsymbol{I}}}_{CB}(k)}{\partial \boldsymbol{p}}= \boldsymbol{C} \frac{\partial \dot{\hat{\boldsymbol{x}}}(k)}{\partial \boldsymbol{p}}+\frac{\partial (\boldsymbol{C}\dot{\boldsymbol{x}})}{\partial \boldsymbol{p}} \bigg\rvert_{\substack{\dot{\boldsymbol{x}}=\dot{\hat{\boldsymbol{x}}}(k)}}, \\
 \\
\frac{\partial \dot{\hat{\boldsymbol{x}}}(k)}{\partial \boldsymbol{p}}= \frac{\partial \Big( \boldsymbol{A}(\boldsymbol{x})\boldsymbol{x} +\boldsymbol{\beta}(\boldsymbol{x}) \Big)}{\partial \boldsymbol{x}}\bigg\rvert_{\substack{\boldsymbol{x}=\hat{\boldsymbol{x}}(k)}} \cdot \frac{\partial \hat{\boldsymbol{x}}(k)}{\partial \boldsymbol{p}} + \frac{\partial \Big( \boldsymbol{A}(\boldsymbol{x})\boldsymbol{x} +\boldsymbol{\beta}(\boldsymbol{x}) \Big)}{\partial \boldsymbol{p}}\bigg\rvert_{\substack{\boldsymbol{x}=\hat{\boldsymbol{x}}(k)}}, \\
 \\
\frac{\partial \hat{\boldsymbol{x}}(k+1)}{\partial \boldsymbol{p}}= \frac{\partial \Big( \big(\boldsymbol{I}+t_s\boldsymbol{A}(\boldsymbol{x})\big)\boldsymbol{x} +t_s \boldsymbol{\beta}(\boldsymbol{x}) \Big)}{\partial \boldsymbol{x}}\bigg\rvert_{\substack{\boldsymbol{x}=\hat{\boldsymbol{x}}(k)}} \cdot \frac{\partial \hat{\boldsymbol{x}}(k)}{\partial \boldsymbol{p}} + \frac{\partial \Big( \big(\boldsymbol{I}+t_s\boldsymbol{A}(\boldsymbol{x})\big)\boldsymbol{x} +t_s \boldsymbol{\beta}(\boldsymbol{x}) \Big)}{\partial \boldsymbol{p}}\bigg\rvert_{\substack{\boldsymbol{x}=\hat{\boldsymbol{x}}(k)}}.
\end{gather*}
For the Input preview case, we can follow a similar procedure, obtaining the same recursive equations detailed above for the forward Euler case, except for the expression of  $\partial \hat{\boldsymbol{x}}(k+1)/\partial \boldsymbol{p}$, which is given by:
$$
\frac{\partial \hat{\boldsymbol{x}}(k+1)}{\partial \boldsymbol{p}}= \frac{\partial \chi (\boldsymbol{x},\boldsymbol{u})}{\partial \boldsymbol{x}}\bigg\rvert_{\substack{\boldsymbol{x}=\hat{\boldsymbol{x}}(k) \\ \tilde{\boldsymbol{u}}(k) \quad \\ \tilde{\boldsymbol{u}}(k+1)}} \cdot \frac{\partial \hat{\boldsymbol{x}}(k)}{\partial \boldsymbol{p}} + \frac{\partial \chi (\boldsymbol{x},\boldsymbol{u})}{\partial \boldsymbol{p}}\bigg\rvert_{\substack{\boldsymbol{x}=\hat{\boldsymbol{x}}(k) \\ \tilde{\boldsymbol{u}}(k) \quad \\ \tilde{\boldsymbol{u}}(k+1)}},
$$
where
$$
\chi (\boldsymbol{x},\boldsymbol{u}) = \left(\boldsymbol{I}-\frac{t_s}{2}\boldsymbol{A}(\boldsymbol{x})\right)^{-1} \left( \Big(\boldsymbol{I}+\frac{t_s}{2}\boldsymbol{A}(\boldsymbol{x})\Big)\boldsymbol{x} + \frac{t_s}{2}\boldsymbol{B}\Big(\boldsymbol{u}(k+1)+\boldsymbol{u}(k)\Big)+\frac{t_s}{2} \Big( \boldsymbol{\beta}(\boldsymbol{x})+\boldsymbol{\beta}(\boldsymbol{x}) \Big) \right).
$$

\bibliographystyle{plain}

\end{document}